\documentclass[11pt,a4paper]{article}
\usepackage{amsmath,amssymb,amsthm,datetime,enumerate,extarrows,fullpage,mathrsfs,mathtools,parskip,tikz}
\usepackage[all,cmtip]{xy}
\usetikzlibrary{arrows,decorations.pathreplacing}

\numberwithin{equation}{section}
\numberwithin{figure}{section}

\begingroup 
\makeatletter 
\@for\theoremstyle:=definition,remark,plain\do{
\expandafter\g@addto@macro\csname th@\theoremstyle\endcsname{
\addtolength\thm@preskip\parskip 
}
}
\makeatother
\endgroup

\newtheorem{theorem}{Theorem}[section]
\newtheorem{lemma}[theorem]{Lemma}
\newtheorem{proposition}[theorem]{Proposition}
\newtheorem{corollary}[theorem]{Corollary}
\theoremstyle{remark}
\newtheorem{remark}[theorem]{Remark}
\theoremstyle{definition}

\newtheorem{example}[theorem]{Example}

\newcommand*{\defeq}{\mathrel{\vcenter{\baselineskip0.5ex \lineskiplimit0pt\hbox{\scriptsize.}\hbox{\scriptsize.}}}=}
\newcommand*{\eqdef}{=\mathrel{\vcenter{\baselineskip0.5ex \lineskiplimit0pt\hbox{\scriptsize.}\hbox{\scriptsize.}}}}

\newcommand{\GL}{\ensuremath{\operatorname{GL}}}
\newcommand{\Tub}{\ensuremath{\operatorname{Tub}}}
\newcommand{\WTub}{\ensuremath{\operatorname{WeakTub}}}
\newcommand{\Emb}{\ensuremath{\operatorname{Emb}}}
\newcommand{\Diff}{\ensuremath{\operatorname{Diff}}}
\newcommand{\Aut}{\ensuremath{\operatorname{Aut}}}
\newcommand{\Map}{\ensuremath{\operatorname{Map}}}
\newcommand{\im}{\ensuremath{\operatorname{im}}}
\newcommand{\id}{\ensuremath{\operatorname{id}}}
\newcommand{\SP}{\operatorname{SP}}

\newcommand{\rlim}{\ensuremath{\underset{{}^{\relbar\joinrel\joinrel\relbar\joinrel\joinrel\rightarrow}}{\lim}}}

\newcommand{\simeqs}{\ensuremath{\stackrel{_s}{\simeq}}}

\newcommand{\longrighthook}{\ensuremath{\lhook\joinrel\relbar\joinrel\rightarrow}}

\renewcommand{\mathbf}[1]{\ensuremath{\boldsymbol{#1}}}

\let\originalleft\left
\let\originalright\right
\renewcommand{\left}{\mathopen{}\mathclose\bgroup\originalleft}
\renewcommand{\right}{\aftergroup\egroup\originalright}

\newdimen{\andlen}
\newdimen{\andskip}
\title{Tubular configurations: equivariant scanning and splitting}
\author{Richard Manthorpe and Ulrike Tillmann}
\date{\monthname\ \the\year}

\begin{document}
	\maketitle



\begin{abstract}
Replacing  configurations of points by configurations of tubular neighbourhoods (or discs)
in a manifold 
$M$ 
we are able to define a natural scanning map that is equivariant 
under the action of the diffeomorphism group of the manifold. 
We also construct
the so-called 
power set map of configuration spaces diffeomorphism equivariantly.
Combining these two constructions yields  stable 
splittings in the sense of Snaith and generalisations thereof 
that are equivariant. In particular
one deduces
stable splittings of homotopy orbit spaces. 
As an application the homology injectivity is proved for diffeomorphism of $M$
that fix an increasing number of points.  
Throughout we work with configurations spaces with 
labels in a fibre bundle over $M$.
\end{abstract}

	\section{Introduction} 
	\label{sec:introduction}


\vskip .1in


There has been much recent interest in configuration spaces of manifolds.  
In one direction,
the work on factorisation algebras and non-commutative Poincare duality by 
Lurie \cite {LUR}, see also Francis \cite {FRA}, is based on the  classical 
work on configuration spaces of May \cite {MILS}, Segal \cite{SILS}, McDuff \cite {MDCS} and Salvatore \cite {MR1851264}. This has
also ignited increased interest in the Goodwillie-Weiss embedding calculus 
\cite {GW}. From this point of view, one is more interested in configurations
of  embedded discs 
than points, and needs to understand the interaction with the 
diffeomorphism group of the background space $M$. 
Our approach to configurations spaces will address both these points.

In another direction, moduli spaces of manifolds and the scanning map have 
been central in the work on the Mumford conjecture and analogues;
see \cite {MWA}, \cite {GAL} and also \cite {TILB} for a survey.
In this context, configuration spaces are moduli spaces of zero dimensional 
manifolds and have provided much intuition.   
The work here was motivated by  some basic question of diffeomorphism
equivariance that arose in this context.

\vskip .1in
{\it Contents and results:}

Let $C_k(M; X)$ denote the space of $k$ unordered, distinct particles in a 
compact smooth manifold $M$ with labels in a pointed space $X$. For a closed 
submanifold $M_0 \subset M$, let $C(M, M_0; X)$ denote the space of 
configurations of particles in $M$ which vanish in $M_0$ 
or at the base-point of 
$X$.

The goal of this paper is two-fold.
First we want to revisit the foundations of the subject and provide a 
natural and equivaraint scanning map that relates the configuration 
spaces  $C(M, M_0;X) $ to  mapping spaces 
or section spaces more generally. The study of these maps
goes back to May \cite{MILS} and Segal \cite {SILS} in the 
case when $M = \mathbb R ^n$, and to McDuff \cite {MDCS} 
and B\"odigheimer \cite {BSS} for general manifolds.
The diffeomorphism group of $M$ acts naturally on both the configuration spaces 
and the section spaces. However, the standard scanning maps, which involve 
choosing a metric on $M$ and are defined by collapsing $\varepsilon$-balls 
around the particles, do not commute with these actions.

Our approach here is to replace a configuration by its space of
tubular neighbourhoods. As the space of tubular neighbourhoods is  
contractible this construction does not 
change the homotopy type of $C_k(M; X)$ and, less obviously, 
also not of $C(M, M_0; X)$. 
This will be proved 
in section 2. In  section 3 the scanning 
map on these enlarged configuration spaces is defined by 
simply collapsing $M$ onto the configuration of 
tubular neighbourhoods. This construction is equivariant under 
the action of the diffeomorphism group; see Theorem 3.8.

\vskip .1in 
Our second goal is to  revisit 
the classical splitting theorems for function spaces
going back to Snaith \cite {SSD} when $M$ is a Euclidean space 
and generalised by B\"odigheimer \cite {BSS} 
to arbitrary manifolds. Using the above results we construct
these splittings equivariantly under the action of the diffeomorphism group; 
the main result in this direction is Theorem 4.8. 
In particular this  gives stable 
splittings of the corresponding homotopy orbit spaces, something that 
for  actions of
compact Lie groups and a restricted class of manifolds 
was previously shown by B\"odigheimer 
and Madsen \cite {BMHQ} by different methods that do not extend to the 
non-compact setting.

The key to the splitting theorem is the construction in section 4.3 of 
diffeomorphism equivariant power set maps
for configuration spaces. This uses the Barratt-Eccles \cite {BEG} 
model for the 
free infinite loop space functor.

In the final section, for connected $M$ with non-empty boundary,
the splitting methods are applied to show that 
the inclusion 
$b: C_k (M; X) \to C_{k+1}(M; X)$, which is well-known to be stably 
split injective, is
indeed equivariantly so. 
As an immediate consequence  we  prove that
a natural homomorphism
$$
\bar b: \Diff (M \smallsetminus \bold k; \partial M) \to \Diff (M 
\smallsetminus \bold {k+1}; \partial M)
$$ 
of diffeomorphisms of $M$ fixing a
set of $ k$ points 
to those fixing a set of  ${k+1}$  points induces a split injection in 
homology on classifying spaces; this is the content of Theorem 4.15.

\vskip .1in
Much of the literature restricts itself to configurations with labels in 
a constant space $X$.  We emphasise that
more generally we consider here  configuration spaces 
with twisted coefficients, 
that is where $X$ is replaced by a fibre bundle $\pi$ over $M$ and
the label space may vary with the points in 
$M$. On the one
hand we will need 
this in the  application we have in mind \cite {TIL} and on the other hand  
it allows us to replace sections spaces with mapping spaces; see Example 4.10.

\vskip .1in
{\it Future work and extensions:}

In
forthcoming work of the second author \cite {TIL}, 
using the results established here, 
the  map $\bar b$ and generalisations thereof will  be 
shown to also  induce  isomorphisms in homology
in a range growing with $k$.

In another direction, the methods of this 
paper can be extended to treat configurations of submanifolds as considered by 
Palmer \cite {PAL} and provide equivariant stable split injections for the 
stable homology isomorphisms in that setting.

	
	\section{Tubular configuration spaces and twisted labels} 
	\label{sec:tubular_configuration_spaces}

We define tubular configuration spaces and show that they are 
homotopic to the usual configuration spaces.
	

	\subsection{The definition of  tubular configuration spaces} 
	\label{sub:tubular_configuration_spaces}
	
	
	Let $M$ be a smooth compact manifold. The configurations 
space of $k$ ordered particles in $M$ is the subspace of 
the $k$-fold Cartesian product of $M$
	\begin{equation*}
		\widetilde C_k(M) \defeq \{(m_1,\ldots,m_k)\in M^k : m_i\neq m_j\, \text{if}\, i\neq j\}.
	\end{equation*}
	Equivalently, $\widetilde C_k(M)$ is the embedding space 
$\Emb(\mathbf{k},M)$, where $\mathbf{k}$ denotes the 0-manifold with $k$ 
points. 
The symmetric group $\Sigma _k$ acts freely and 
the configuration space of $k$ unordered particles in $M$ is the 
orbit space $C_k(M)$. 
When $k=0$ there is only one configuration,  the empty configuration, and
 $\widetilde C_0(M) = C_0(M) = *$.
	
	Let $M_0\subset M$ be a (possibly empty) compact submanifold. The configuration space of particles in $M$ modulo $M_0$ is then defined as
	\begin{equation*}
		C(M,M_0)\defeq \left(\coprod_{k = 0}^{\infty} C_k(M)\middle)\!\middle/\! \sim \right.
	\end{equation*}
	where $(m_1,\ldots,m_k)\sim (m_1,\ldots,m_{k-1})$ if $m_k\in M_0$. We think of this relation as particles vanishing in $M_0$. 

To define the tubular configuration spaces of particles in $M$ 
we replace a configuration $\mathbf{m} = (m_1,\ldots,m_k)$ by the space of tubular neighbourhoods of its particles considered as a 0-dimensional submanifold of $M$.
	
	Let $W$ be a manifold without boundary containing $M$ 
as a codimension zero submanifold. More precisely, 
if $M$  has empty boundary let $W = M$ and  
otherwise let $W = M \cup (\partial M\times [0,1))$ 
be $M$ with an open collar attached. Similarly we define $M^+$ as
$M^+ = M = W$ when $M$ has no boundary and as $M^+= M\cup (\partial M
\times [0, 1/2))$ otherwise.
Such  manifolds $M^+$  and $W$ are required so that particles on the 
boundary of $M$ admit tubular neighbourhoods. 
	
	Let $P\subset M$ be a neat submanifold, let $\nu$ be
its normal bundle and identify $P$ with the image of the zero section in $\nu$. By a tubular neighbourhood of $P$ in $M$ we mean an embedding $f:\nu\rightarrow M$ which restricts to the identity on the zero-section and for which the composition
	\begin{equation}
		\tag{\id}
		\nu\longrighthook \nu\oplus TP = T\nu|_P\overset{df}{\longrightarrow}TM|_P\longrightarrow \nu
	\end{equation}
	is the identity on $\nu$. We call this last property $(\id)$. Denote the space of tubular neighbourhoods of $P$ by $\Tub(P)$ and topologise it as a subspace of the embedding space $\Emb(\nu,M)$ with the $C^\infty$ topology.
	
	For an ordered or unordered configuration 
$\mathbf{m}$ of $M$ let $\Tub(\mathbf{m})$ denote the space of 
tubular neighbourhoods of $\mathbf{m}$ considered as a 0-dimensional 
submanifold of $M^+$. 
Define the tubular configuration space of $k$ ordered particles in $M$ as the disjoint union
	\begin{equation*}
		\widetilde E_k(M) \defeq \coprod_{\mathclap{\mathbf{m}\in \widetilde C_k(M)}} \Tub(\mathbf{m})
	\end{equation*}
	equipped with an appropriate topology which restricts to the 
ordinary $C^\infty$ topology on the fibres. 
We defer a description of the topology to Section 
\ref{sub:the_parc_C_topology}. 
We refer to a tubular neighbourhood $f\in\widetilde E_k(M)$ as a 
tubular configuration and we often write it as a collection of its 
components $(f_1,\ldots,f_k)$, with $f_i:T_{m_i}M\rightarrow M^+$. 
The symmetric group on $k$ points acts freely on 
$\widetilde E_k(M)$ and we define the tubular configuration space 
of $k$ unordered particles as the orbit space $E_k(M)$.

	Analogous to the ordinary configuration spaces, we define the tubular configuration space of particles in $M$ modulo $M_0$ as
	\begin{equation*}
		E(M,M_0)\defeq \left(\coprod_{k = 0}^{\infty} E_k(M)\middle)\!\middle/\! \sim \right.
	\end{equation*}
	where $(f_1,\ldots,f_k)\sim (f_1,\ldots,f_{k-1})$ if $f_k(0)\in M_0$. Components of a tubular configuration vanish if their midpoint is in $M_0$.

	
 \begin{figure}[htb]
                \begin{center}
                        \begin{tikzpicture}
                                \draw (0,0) node
{\includegraphics{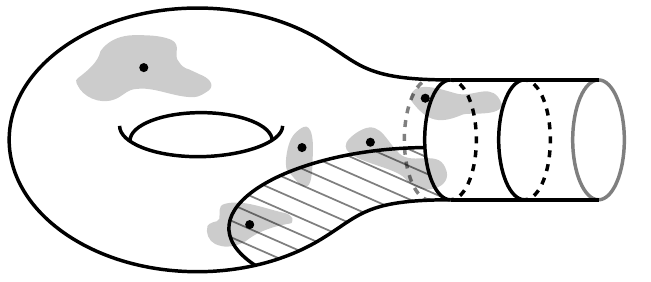}};
                                \draw (0.2,-0.5) node {$M_0$};
                                \draw[decorate,
decoration={brace,amplitude=4pt}] (-3.21,1.5) -- (1.27,1.5)
node[align=center, midway, anchor=south, yshift=3pt] {$M$};
                                \draw[decorate,
decoration={brace,amplitude=4pt}] (1.28,1.5) -- (2.78,1.5)
node[align=center, midway, anchor=south, yshift=3pt] {$W\smallsetminus M$};
                                \draw[decorate,
decoration={brace,amplitude=4pt}] (2.02,-1.5) -- (-3.21,-1.5)
node[align=center, midway, anchor=north, yshift=-3pt] {$M^+$};
                        \end{tikzpicture}
                \end{center}
                \label{fig:tubular_configuration}
                \caption{\it The image of a tubular configuration of five
particles on a punctured torus. The underlying particles are marked as black
dots.}
        \end{figure}
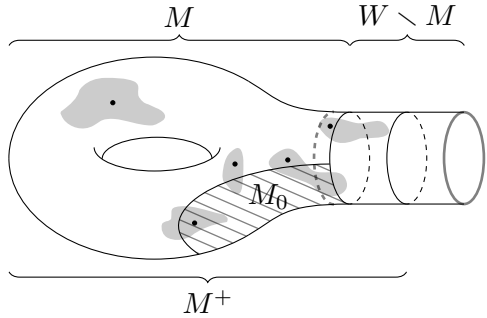

\vskip .1in
	
	\begin{remark}
		As an alternative approach, consider configuration spaces with neighbourhoods of each particle without an identification with the normal bundle. We think of this as configurations of submanifolds diffeomorphic to a finite union of open disks with one marked point. Note that for a configuration $\mathbf{m}$ of $k$ particles in $M$, there are isomorphisms
		\begin{gather*}
			\Tub(\mathbf{m}) \cong \WTub(\mathbf{m})\big/\textstyle{\bigoplus_i}\GL_n(\mathbb{R})\\
			\text{and}\\
			\Tub(\mathbf{m})\big/\textstyle{\bigoplus_i}\Diff(T_{m_i}M,\id)\
				\cong \WTub(\mathbf{m})\big/\textstyle{\bigoplus_i}\Diff(T_{m_i}M,0)
		\end{gather*}
		where $\WTub(\mathbf{m})$ is the space of weak 
tubular neighbourhoods, that is the space of embeddings of the normal bundle of $\mathbf{m}$ into $M^+$ which 
restrict to the identity on the zero section but do not satisfy property $(\id)$, $\Diff(T_{m_i}M,0)$ is the group of diffeomorphisms which fix $0$ and $\Diff(T_{m_i}M,\id)\subset\Diff(T_{m_i}M,0)$ are the diffeomorphisms whose derivatives at zero are the identity. 
The second space is the space of embedded discs with given midpoints.
We note that 
$$
\Diff ( T_{m_i}M,0) \simeq \GL _n(\mathbb R) \text { and } 
\Diff ( T_{m_i} M, \id ) \simeq *.
$$
Hence, the space of embedded discs (with fixed midpoint) is homotopic to the
space of tubular neighbourhoods, in other words the two spaces 
above are homotopic.
	\end{remark}
	

	\subsection{The parc $C^\infty$ topology} 
	\label{sub:the_parc_C_topology}

	We generalise the compact-open topology for spaces of partial maps with closed domain as defined in \cite{BBPM} to a $C^\infty$ topology for spaces of smooth partial maps. This generalisation allows us to topologise tubular configuration spaces as subspaces of certain smooth partial mapping spaces.

	Let $X$ and $Y$ be topological spaces. A partial map $f:X\rightarrow Y$ is a map $A\rightarrow Y$ for some subspace $A\subseteq X$. We call $A$ the domain of $f$ and denote it by $\mathcal{D}(f)$. A parc map or partial map with closed domain is a partial map $f$ such that $\mathcal{D}(f)$ is closed in $X$. Let $P_c(X,Y)$ denote the set of parc maps $X\rightarrow Y$. The parc mapping space is $P_c(X,Y)$ equipped with the parc compact-open topology which has sub-base the sets
	\begin{equation*}
		(K,U)\defeq \{f\in P_c(X,Y) : f(K\cap \mathcal{D}(f))\subseteq U\}
	\end{equation*}
	for all compact sets $K\subseteq X$ and open sets $U\subseteq Y$.
	
	Similarly there is a paro mapping space $P_o(X,Y)$ for partial maps with open domain \cite{ABPM}. $P_o(X,Y)$ is equipped with the paro compact-open topology with sub-base the sets
	\begin{equation*}
		(K,U)\defeq \{f\in P_o(X,Y) : K\subseteq \mathcal{D}(f),\; f(K)\subseteq U\}.
	\end{equation*}
	
	\begin{remark}
		Given a fixed closed (open) set $A\in X$, the topology on the subset of parc (paro) maps $X\rightarrow Y$ which are defined precisely on $A$ coincides with the ordinary compact-open topology on the mapping space $\Map(A,Y)$.
	\end{remark}
	
	Let $M$ and $N$ be $C^r$ manifolds with $r<\infty$ and denote the set of partial $C^r$ maps $M\rightarrow N$ with closed domain by $C^rP_c(M,N)$. We equip this set with a generalisation of the $C^r$ topology as follows. 
Let $f\in C^rP_c(M,N)$ be a parc $C^r$ map and let $(\phi,U)$ and $(\psi,V)$ be charts for $M$ and $N$ respectively. Let $K\subseteq U$ be a compact set such that $f(K\cap \mathcal{D}(f))\subseteq V$ and let $0<\varepsilon\le\infty$. Define a parc sub-basic neighbourhood
	\begin{equation*}
		\mathcal{N}^r(f;(\phi,U),(\psi,V),K,\varepsilon)
	\end{equation*}
	to be the set of parc $C^r$ maps $g:M\rightarrow N$ such that $g(K\cap \mathcal{D}(g))\subseteq V$ and
	\begin{equation*}
		\|D^k(\psi f\phi^{-1})(x) - D^k(\psi g\phi^{-1})(x)\|<\varepsilon
	\end{equation*}
	for all $x\in \phi(K)$ and $k = 0,\ldots, r$. The set of all such neighbourhoods form a sub-base for the parc $C^r$ topology on $C^rP_c(M,N)$. For $C^\infty$ manifolds $M$ and $N$, $C^\infty P_c(M,N)$ is the space of parc $C^\infty$ maps equipped with the parc $C^\infty$ topology. This is simply the union of the topologies induced by the inclusions $C^\infty P_c(M,N)\hookrightarrow C^rP_c(M,N)$ for $r$ finite.

	\begin{remark}
		\label{rmk:parc_C_top_agrees}
		It follows immediately from the definition that the subspace of $C^\infty P_c(M,N)$ consisting of maps defined on a smooth closed submanifold $A$ is the smooth mapping space $C^\infty(A,N)$ equipped with the ordinary $C^\infty$ topology.
	\end{remark}
	
	\begin{remark}
		The parc smooth mapping space is functorial in each argument. Let $M$, $N$ and $Q$ be smooth manifolds, then smooth maps $\phi:N\rightarrow Q$ and $\psi:M\rightarrow Q$ induce the following continuous maps.
		\begin{equation*}
			\begin{aligned}
				\phi_*: C^\infty P_c(M,N) &\longrightarrow C^\infty P_c(M,Q) \\
				f &\longmapsto \phi\circ f
			\end{aligned}
			\hspace{0.1\textwidth}
			\begin{aligned}
				\psi^*: C^\infty P_c(M,N) &\longrightarrow C^\infty P_c(Q,N) \\
				f &\longmapsto f\circ\psi|_{\psi^{-1}(\mathcal{D}(f))}
			\end{aligned}
		\end{equation*}
	\end{remark}
	
	
	\subsection{The topology of tubular configuration spaces} 
	\label{sub:topologising_E}
	
	Returning to configuration spaces, let $M \subset M^+ \subset W$ 
be as defined in section 2.1, and define the tubular configuration space of $k$ unordered particles in $M$ as the subspace
	\begin{equation*}
		E_k(M)\defeq\{f\in C^\infty P_c(TM,M^+) : \mathcal{D}(f) = \coprod_iT_{m_i}M \text{ and } f\in\Tub(\mathbf{m}) \text{ for some } \mathbf{m}\in C_k(M)\}.
	\end{equation*}
	This  agrees set-wise with the previous definition. 
Moreover, for any $\mathbf{m} \in C_k(M)$ the topology on the subspace 
$\Tub(\mathbf{m})\subset E_k(M)$ is compatible 
with the ordinary $C^\infty$ topology by Remark \ref{rmk:parc_C_top_agrees}.
	

	\begin{lemma}
		The projection $p:E_k(M)\rightarrow C_k(M)$, $f\mapsto (f_1(0),\ldots,f_k(0))$ is continuous.
	\end{lemma}
	
	\begin{proof}
		Let $r:W\rightarrow M$ be the retract collapsing the collar onto $\partial M$ and let $j:M\hookrightarrow TM$ be the zero section. Note that $C_k(M)$ can be identified with the space of partial (smooth) identity maps $M\rightarrow M$ with domain a finite subset of size $k$. 
Then $p$ is the restriction of the induced map 
$$
C^\infty P_c (TM, M^+) 
\xrightarrow{(r|_{M^+})_*}C^\infty P_c(TM,M)
\xrightarrow{j^*} C^\infty P_c(M,M)$$ 
to $E_k(M)$, 
which is 
continuous by Remark 2.4.
	\end{proof}
	
	\begin{proposition}
		\label{prop:E_C_bundle}
		$p:E_k(M)\rightarrow C_k(M)$ is a fibre bundle.
	\end{proposition}
	
	\begin{proof}
		Let $\mathbf{m} = (m_1,\ldots,m_k)$ be a configuration 
in $C_k(M)$. 
We construct a local trivialisation for $p$ around $\mathbf{m}$. 
Let $h:TM|_{\mathbf{m}}\rightarrow M^+$ be a tubular 
neighbourhood of $\mathbf{m}$ and let $DM$ and 
$\mathring DM$ be the closed and open disk subbundles of 
$TM$ for some fixed metric.
Then the restriction of $h$ to $DM|_{\mathbf{m}}$ is a closed tubular neighbourhood. Choose a continuous family of diffeomorphisms of $h(DM|_{\mathbf{m}})$ fixing the boundary
		\begin{equation*}
			\tau:\Gamma(\mathring DM|_{\mathbf{m}}) \longrightarrow \Diff(h(DM|_{\mathbf{m}}),\partial)
		\end{equation*}
		parameterised by the space of sections of $\mathring DM|_{\mathbf{m}}$ and such that for each section $s\in \Gamma(\mathring DM|_{\mathbf{m}})$, $\tau_s\circ h\circ s = \id_{\mathbf{m}}$ where $\tau_s\defeq \tau(s)$. 
We assume the diffeomorphisms are extended to $M^+$ and $W$ 
by fixing $W\smallsetminus h(DM|_{\mathbf{m}})$.
		
		There is a homeomorphism 
		\begin{align*}
			\phi:\Gamma(\mathring DM|_{\mathbf{m}}) &\longrightarrow \left(\prod_ih(D_{m_i}M)\middle)\!\middle/\!\Sigma_k\right. \eqdef V\subseteq C_k(M) \\
			s &\longmapsto (h\circ s(m_1),\ldots,h\circ s(m_k))
		\end{align*}
		whose image is an open neighbourhood of $\mathbf{m}$ in $C_k(M)$. Intuitively, $h$ defines an open ball around each particle in $\mathbf{m}$ and $V$ is the set of all configurations with precisely one particle in each open ball. For each configuration $\mathbf{n}\in V$, $\sigma_{\mathbf{n}}\defeq\tau\circ \phi^{-1}(\mathbf{n})$ is the diffeomorphism of $W$ which moves the particles in $\mathbf{n}$ onto the particles in $\mathbf{m}$.
		
		Local trivialisations over $V$ are then given by
		\begin{align*}
			E_k(M)|_V &\longmapsto V\times\Tub(\mathbf{m}) \\
			\Tub(\mathbf{n})\ni f &\longmapsto \left(\mathbf{n}, \sigma_{\mathbf{n}}\circ f\circ \left(d\sigma_{\mathbf{n}}\right)^{-1}\right).
		\end{align*}
		Here $\left(d\sigma_{\mathbf{n}}\right)^{-1}$ maps the normal bundle of $\mathbf{m}$ to the normal bundle of $\mathbf{n}$, $f$ maps this to a neighbourhood of $\mathbf{n}$ and $\sigma_{\mathbf{n}}$ maps this neighbourhood to a neighbourhood of $\mathbf{m}$, thus the composition is indeed a tubular neighbourhood of $\mathbf{m}$.
	\end{proof}

		For each $k\ge 0$ there is a covering map $q:\widetilde C_k(M)\rightarrow C_k(M)$. We define the tubular configuration space of $k$ ordered particles in $M$ as the pullback $\widetilde E_k(M)\defeq q^*E_k(M)$ and let $\tilde p: \widetilde E_k (M) \to \widetilde C_k (M)$
be the corresponding fibre bundle.
	
The projection $p$ has a right inverse $\sigma$. To define $\sigma$	
	choose  a Riemannian metric  on $W$. 
For a configuration $\mathbf{m}$ let $\varepsilon_1$ be the smallest 
distance in $W$ between any two of its particles and let $\varepsilon_2$ 
be the greatest value such that for each $i$, 
$\exp_{m_i}:B_{\varepsilon_2}(T_{m_i}M,0)\rightarrow B_{\varepsilon_2}
(M^+,m_i)$ is a diffeomorphism. 
Let $\varepsilon_{\mathbf{m}} = \min\{\varepsilon_1,\varepsilon_2\}$, 
then $\mathbf{m}\mapsto \varepsilon_{\mathbf{m}}$ is a continuous map 
$\widetilde C_k(M) \rightarrow \mathbb{R}_{>0}$. 
For each $i$ define $f_i:T_{m_i}M\rightarrow M^+$ 
by $v\mapsto \exp_{m_i}\left(\frac{2\varepsilon_{\mathbf{m}}v}
{\pi|v|}\arctan|v|\right)$, then $(m_1,\ldots,m_k)\mapsto (f_1,\ldots,f_k)$ is a section $\tilde \sigma:\widetilde C_k(M)\rightarrow \widetilde E_k(M)$. Moreover, this section is $\Sigma_k$-equivariant and descends to a section $\sigma:C_k(M)\rightarrow E_k(M)$.

	\begin{corollary}
		The projections $p$ and $\tilde p$ are homotopy equivalences with homotopy inverses given by the global sections $\sigma$ and $\tilde \sigma$.
	\end{corollary}
	
	\begin{proof}
		We have seen that the projections $p$ and $\tilde p$ 
give tubular configuration spaces the structure of fibre bundles. 
The fibre over any configuration $\mathbf{m}$ is the space of 
tubular neighbourhoods $\Tub(\mathbf{m})$. The space of tubular 
neighbourhoods of any compact submanifold is contractible, see for example
 \cite{GHST}.  The contractions of each fibre determine homotopies $\sigma \circ p\sim\id_{E_k(M)}$ and $\tilde \sigma\circ\tilde p\sim\id_{\widetilde E_k(M)}$.
	\end{proof} 

	
	\subsection{Twisted labels and homotopy equivalences} 
	\label{sub:twisted_labels}
	
	It is common to add local data to configurations in the form of labels in a parameter space. In this paper we consider an extension of this notion which allows the parameter space to vary as the fibre of a fibre bundle over the underlying manifold. We say the configurations have \emph{twisted labels}. 
	
	Let $M$ be a smooth compact manifold and let $\pi:Y\rightarrow M$ be a fibre bundle and a zero section $o:M\rightarrow Y$. 
Furthermore, assume that for each $m\in M$, the fibre $Y_m$ over $m$
 is well-pointed with base-point $o(m)$. We define the configuration space of $k$ ordered particles in $M$ with twisted labels in $\pi$ as
	\begin{equation*}
		\widetilde C_k(M;\pi)\defeq \{(\mathbf{m},x)\in \widetilde C_k(M)\times P_c(M,Y):\mathcal{D}(x) = \mathbf{m},\, x(m_i)\in Y_{m_i}\}
	\end{equation*}
or equivalently 
$\widetilde C_k(M;\pi)\defeq \{(\mathbf{m};\mathbf{x})\in \widetilde C(M)\times Y^k : \pi(x_k) = m_k\}$. We prefer here the definition in terms of 
partial maps as it motivates the definition of tubular configuration spaces 
with twisted labels below.
	
	\begin{example}
		An important example is when $Y$ is the trivial bundle $M\times X$ where $X$ is well-pointed with base-point $*$ and the zero section is $o(m) = (m,*)$. In this case we write $\widetilde C_k(M;X)$ for 
 the configuration space with labels in $X$.
	\end{example}
	
	As in the unlabelled case, we define $C_k(M;\pi)$, the configuration space of $k$ unordered particles with twisted labels in $\pi$, 
as the orbit space under the obvious $\Sigma_k$-action and $C_0(M;\pi) = *$. Now let $M_0\subset M$ be a compact submanifold, then the configuration space of particles in $M$ modulo $M_0$ with twisted labels in $\pi$ is
	\begin{equation*}
		C(M,M_0;\pi)\defeq \left(\coprod_{k = 0}^{\infty} C_k(M;\pi)\middle)\!\middle/\! \sim \right.
	\end{equation*}
	where $(m_1,\ldots,m_k;x_1,\ldots,x_k)\sim (m_1,\ldots,m_{k-1};x_1,\ldots,x_{k-1})$ if $m_k\in M_0$ or $x_k= o|_{m_k}$. Here $x_i$ denotes the $i$th component of $x$, i.e. $x_i = x|_{m_i}$. 
	
	Let $M \subset M^+ \subset W$ be as in section \ref{sub:tubular_configuration_spaces} and let $(Y,\pi,o)$ 
be extended over $W$ as the pullback along the retract $r: W\rightarrow M$ 
which maps the collar to $\partial M$. The tubular configuration space of $k$ ordered particles in $M$ with twisted labels in $\pi$ is
	\begin{equation*}
		\widetilde E_k(M;\pi)\defeq \{(f,s)\in \widetilde E_k(M)\times P_o(M^+,Y):\mathcal{D}(s) = \im(f),\, \pi s = \id_{\im(f)}\}
	\end{equation*}
	and the unordered tubular configuration space is the orbit space $E_k(M;\pi)$.
	
	The following  propositions are proved using similar 
techniques as in the proof of Proposition \ref{prop:E_C_bundle} 
and Corollary 2.7.
	
	\begin{proposition}
		The projection $q:E_k(M;\pi)\rightarrow E_k(M)$, $(f;s)\mapsto f$ is a fibre bundle with fibre $q^{-1}(f) = \Gamma(Y|_{\im(f)})$, 
the space of sections of $Y$ over the image of $f$. \qed
	\end{proposition}
	
	\begin{proposition}
		The projection $p:E_k(M;\pi)\rightarrow C_k(M;\pi)$, $(f;s)\mapsto (p(f);s|_{p(f)})$ is a fibre bundle with contractible fibres. \qed
	\end{proposition}

	Analogous to the ordinary configuration spaces, we define the tubular configurations space of particles in $M$ modulo $M_0$ with twisted labels in $Y$ 
to be
	\begin{equation*}
		E(M,M_0;\pi)\defeq \left(\coprod_{k = 0}^{\infty} E_k(M;\pi)\middle)\!\middle/\! \approx \right.
	\end{equation*}
	where $(f_1,\ldots,f_k;s_1,\ldots,s_k) 
\approx (f_1,\ldots,f_{k-1};s_1,\ldots,s_{k-1})$ 
if $f_k(0)\in M_0$ or $s_k = o|_{\im(f_k)}$.
We emphasise that we require the entire section $s_k$ to agree 
with the zero section rather than just $s_k (m_k) = o(m_k)$. 
Our choice of definition here is motivated by 
the construction of the scanning map in the next section.
Note that 
whereas a particle in a configuration vanishes if its label agrees with the 
zero section, a component of a tubular configuration vanishes only if the entire section over the image of that component agrees with the zero section. A 
tubular configuration does not necessarily vanish 
when the underlying particle in the configuration vanishes. 
This fact complicates the proof of the next result. 
	
	\begin{proposition}
		\label{prop:proj_weak_equiv}
		There is a well defined weak homotopy equivalence $p:E(M,M_0;\pi)\rightarrow C(M,M_0;\pi)$ induced by the projections $p:E_k(M;\pi)\rightarrow C_k(M;\pi)$.
	\end{proposition}
	
	
	\begin{proof}
		Define $p:E(M,M_0;\pi)\rightarrow C(M,M_0;\pi)$ by $p[(f;s)]\defeq [p(f;s)]$ for any representative of the class. This is well defined since if a component of a configuration in $E(M,M_0;\pi)$ vanishes, the underlying particle in $C(M,M_0;\pi)$ does. The number of particles in configurations and tubular configurations induce filtrations
		\begin{align*}
			*= C^0(M,M_0;\pi)\subset C^1(M,M_0;\pi)\subset\ldots\subset C(M,M_0;\pi),
			\quad
			C^n(M,M_0;\pi)\defeq \left(\coprod_{k = 0}^{n} C_k(M;\pi)\middle)\!\middle/\! \sim \right. \phantom{.}\\
			\shortintertext{\hspace{\andskip}and}
			*= E^0(M,M_0;\pi)\subset E^1(M,M_0;\pi)\subset\ldots\subset E(M,M_0;\pi),
			\quad
			E^n(M,M_0;\pi)\defeq \left(\coprod_{k = 0}^{n} E_k(M;\pi)\middle)\!\middle/\! \approx \right.
		\end{align*}
		and $p$ respects the filtrations. Topologically the configuration spaces $C(M,M_0;\pi)$ and $E(M,M_0;\pi)$ are the colimits of the filtration sequences. For each $k$ define the subspaces of configurations in $C_k(M;\pi)$ and $E_k(M;\pi)$ in which at least one particle or component vanishes under the equivalence relation $\sim$ or $\approx$ respectively. More precisely, define
		\begin{gather*}
			B_{C_k}(M,M_0;\pi) \defeq\{(\mathbf{m},\mathbf{x}) 
\in C_k (M; \pi) :m_k\in M_0 \text{ or } x_k = o|_{m_k}\}\\
			\text{and}\\
			B_{E_k}(M,M_0;\pi) \defeq\{(f,s)
\in E_k (M; \pi):f_k(0)\in M_0 \text{ or } s_k = o|_{\im(f_k)}\}.
		\end{gather*}
		For each $k\ge 1$ we have a commutative diagram
		\begin{equation*}
			\xymatrix@!0@R=25pt@C=70pt{
			 & E_k(M;\pi) \ar'[d][dd] \ar[rr]
			   & & C_k(M;\pi) \ar[dd]
			\\
			 B_{E_k}(M,M_0;\pi) \ar[ur]\ar[dd] \ar[rr]
			 & & B_{C_k}(M,M_0;\pi) \ar[ur]\ar[dd]
			\\
			 & E^k(M,M_0;\pi) \ar'[r][rr]
			   & & C^k(M,M_0;\pi)
			\\
			E^{k-1}(M,M_0;\pi) \ar[ur] \ar[rr]
			 & & C^{k-1}(M,M_0;\pi) \ar[ur]
			}
		\end{equation*}
		in which all the maps from left to right are induced by $p$, 
the  vertical maps send configurations to their equivalence classes, and 
the maps from front to  back  are inclusions. 
The upper maps from front to back are cofibrations by our condition 
on the section $o$, and 
the left and right faces are homotopy pushout squares. Thus if $p$ 
induces weak homotopy equivalences on the two upper maps and the front map, 
it also will induce a weak homotopy equivalence on the lower back map. 
		
		By the previous proposition
the projections $p:E_k(M;\pi)\rightarrow C_k(M;\pi)$ 
are weak 
homotopy equivalences for each $k$. 
Note that $E^0(M,M_0;\pi) = * = C^0(M,M_0;\pi)$. Thus, if $p$ induces weak 
homotopy equivalences $B_{E_k}(M,M_0;\pi)\rightarrow B_{C_k}(M,M_0;\pi)$ 
for each $k$ we can proceed by induction on $k$ to show that
$E^k(M,M_0;\pi)\rightarrow C^k(M,M_0;\pi)$ is a weak homotopy equivalence and
hence obtain the result.

                To see  that $B_{E_k}(M,M_0;\pi)\rightarrow B_{C_k}(M,M_0;\pi)$
is a  homotopy equivalence,  define subspaces
                \begin{gather*}
                        B^{o}_{C_k} \defeq \{(\mathbf{m};\mathbf{x}):x_k=o|_{m_k}\},\,
                        B^{M_0}_{C_k} \defeq \{(\mathbf{m};\mathbf{x}):m_k\in M_0\} \subset B_{C_k}(M,M_0;\pi) \\
                        \text{and}\\
                        B^{o}_{E_k}  \defeq \{(f;s):s_k = o|_{\im(f)}\},\, 
                        B^{M_0}_{E_k} \defeq \{(f;s):f_k(0)\in M_0\} \subset B_{E_k}(M,M_0;\pi).
                \end{gather*}
                There is a commutative diagram
                \begin{equation*}
                        \xymatrix@!0@R=25pt@C=60pt{
                         & B^o_{E_k} \ar'[d][dd] \ar[rr]
                           & & B^o_{C_k} \ar[dd]
                        \\
                         B^o_{E_k}\cap B^{M_0}_{E_k} \ar[ur]\ar[dd] \ar[rr]
                         & & B^o_{C_k}\cap B^{M_0}_{C_k} \ar[ur]\ar[dd]
                        \\
                         & B_{E_k}(M,M_0;\pi) \ar'[r][rr]
                           & & B_{C_k}(M,M_0;\pi)
                        \\
                        B^{M_0}_{E_k} \ar[ur] \ar[rr]
                         & & B^{M_0}_{C_k} \ar[ur]
                        }
                \end{equation*}
                where the maps from left to right are induced by $p$ and the left and right faces are homotopy pushout squares. Note that $B^{o}_{E_k}\subset E_k(M;\pi)|_{B^o_{C_k}}$ are subbundles of $E_k(M;\pi)\rightarrow E_k(M)$ with fibres $\Gamma(X|_{\im(f)\smallsetminus \im(f_k)})$ and $\{s\in\Gamma(X|_{\im(f)}):s_k|_{f_k(0)} = o|_{f_k(0)}\}$ respectively over $f$. The fibres are homotopy equivalent so $B^o_{E_k}\rightarrow B^o_{C_k}$ is a weak homotopy equivalence. The front maps from left to right are weak homotopy equivalences since $B^{M_0}_{E_k} = E_k(M;\pi)|_{B^{M_0}_{C_k}}$ and $B^o_{E_k}\cap B^{M_0}_{E_k} = B^o_{E_k}|_{B^{M_0}_{C_k}}$ as bundles over $C_k(M;\pi)$. So $B_{E_k}(M,M_0;\pi)\rightarrow B_{C_k}(M,M_0;\pi)$ is a weak homotopy equivalence.
	\end{proof}

\begin{remark}
If $Y$ has the homotopy type of a CW-complex then the spaces $C_k(M; \pi)$ and
$C(M, M_0; \pi) $, and their filtration subspaces and quotients  
have the homotopy types of
CW-complexes \cite {MDCS} 
and hence the same holds for  tubular configuration spaces. 
The map in the previous proposition is thus in fact a homotopy equivalence.
\end{remark}

	
	
	\section{The scanning map} 
	\label{sec:the_scanning_map}

We define the scanning map for tubular configuration spaces, show that it 
induces a homotopy equivalence under the usual assumptions 
and that it is equivariant with respect to the group of those
diffeomorphisms of $M$ 
which can be lifted to an isomorphism of the bundle $\pi: Y \to M$.

	
	\subsection{Definition and homotopy equivalence} 
	\label{sub:definition}
	Let $M_0 \subset M \subset M^+ \subset W$ and $(Y,\pi,o)$ be as defined in the 
previous section, and let $n=\dim(M)$. Denote by $\dot TW$ 
the fibrewise one point compactification of the tangent bundle of $W$ 
and let $\dot T_zW$ denote the fibre of $\dot TW$ over a point $z \in W$. 
Let $\dot\tau_\pi \defeq \dot TW\wedge_W Y$ be the fibrewise smash product of $\dot TW$ with the bundle $Y$, 
where the base-point in each fibre of $Y$ is determined by the zero-section $o$.
Then, if $X$ denotes a typical fibre of $\pi$, 
$\dot\tau_\pi$ is a $\Sigma^n X$-bundle over $W$ 
with a canonical section $*:W\rightarrow \dot\tau_\pi$. 
Let $\Gamma(W\smallsetminus M_0,W\smallsetminus M;\pi)$ denote the 
space of sections of $\dot\tau_\pi$ which are defined outside of $M_0$ 
and which agree with $*$ outside of $M$ (and hence on $\partial M$). In this section we define a scanning map $\gamma:E(M,M_0;\pi)\rightarrow \Gamma(W\smallsetminus M_0,W\smallsetminus M;\pi)$.
	
	We begin by constructing a sequence of maps $\gamma_k:E_k(M;\pi)\rightarrow \Gamma(W\smallsetminus M_0,W\smallsetminus M;\pi)$. 
Intuitively, our scanning map will send a tubular configuration $(f;s)\in E_k(M;\pi)$ to a section defined by $s$ on the $Y$ component 
and by a modification  of $f^{-1}$ over $W$ on the $\dot TW$ component, sending points outside of $\im(f)$ to the compactification points in the appropriate fibres.
	
	\begin{lemma}
		\label{lem:canonical_iso}
\cite{LSM}
		Let $V$ be a finite dimensional vector space. Given a vector $v\in V$, there is a canonical linear isomorphism $V\rightarrow T_vV$. Moreover, for any finite dimensional vector space $U$ and any linear map $L:V\rightarrow U$ the following diagram commutes.
		\begin{equation*}
			\xymatrix{
			V \ar[r]^\cong \ar[d]_L & T_vV \ar[d]^{d_vL} \\
			W \ar[r]_\cong & T_{Lv}W
			}
		\end{equation*}
\qed
	\end{lemma}
	
	For any $(f,s)\in E_k(M;\pi)$ and $z\in \im(f)$, $f^{-1}(z)\in T_{m_i}M$ for some $i$. For each $z$ and each $f$ let $\phi_{f,z}:T_{m_i}M\rightarrow T_{f^{-1}(z)}(T_{m_i}M)$ be the canonical isomorphism. 
$\phi_{f,z}$ varies continuously in $f$ and $z$. 
Let $(f,s) \in E_k (M; \pi)$. Having removed any components 
$(f_i, s_i)$ with $m_i \in M_0$ if necessary, we define a section of $\dot\tau_\pi$ by
	\begin{equation*}
		z\mapsto
		\begin{cases}
			\left(d_{f^{-1}(z)}f\circ\phi_{f,z}\circ f^{-1}(z),s(z)\right) & \text{if } z\in\im(f) \cap ( M^+\smallsetminus M_0)\\
			*_z &\text{otherwise.}
		\end{cases}
	\end{equation*}

Although the above construction gives a well-defined map, it is not continuous.
The problem is that the contribution of $f_i$ in the section  of 
$\dot\tau_\pi$ suddenly 
vanishes as $m_i$ reaches $M_0$. We introduce the following 
modification of the above map for components $(f_i, s_i)$ with $m_i$ in a
collar neighbourhood of $\partial M_0$ in $M$ which makes sure that for a fixed 
 $z \in 
\im (f_i)$ its image under the section
goes to the point at infinity in $\dot T_zW$  
as $m_i$ gets close 
to the boundary of $M_0$ in $M$. To this end we
fix a collar and choose a metric on $W$ such that the chosen 
collar has width 1. Let $\delta (m_i)$
be the distance of $m_i$ to $M_0$. We now multiply the above formula by  
$\exp (\frac 1 {\delta (m_i)}) $ when $z \in \im (f_i) $ and $m_i $ is in the collar. 
Note that we only use 
the metric on the collar. The resulting map $\gamma ^+$ is now continuous. 
  
Note that the target of $\gamma^+$ is a section that may not vanish outside $M$
as our tubular discs may have image in the larger $M^+$. Let $\rho:
\Gamma (W\smallsetminus M_0, W\smallsetminus M^+;\pi) \to \Gamma (W\smallsetminus M_0, W\smallsetminus M;\pi)$ 
be the homotopy equivalence induced by 
a diffeomorphism of $W$ that maps $\partial M \cup [0,1/2)$ into a collar
of $\partial M$ in $M$ and leaves $M$ away from the collar pointwise fixed.
A homotopy inverse is given by the inclusion
$\Gamma (W\smallsetminus M_0, W\smallsetminus M;\pi) \to \Gamma (W\smallsetminus M_0, W\smallsetminus M^+;\pi)$ 
Define
$$
\gamma:= \rho \circ \gamma ^+ : E(M,M_0;\pi)\rightarrow \Gamma(W\smallsetminus M_0,W\smallsetminus M;\pi).
$$

	\begin{remark}
		\label{rmk:scan_homotopic}
		The composition
		\begin{equation*}
			C(M,M_0;\pi) \xlongrightarrow{\sigma} E(M,M_0;\pi) \xlongrightarrow{\gamma} \Gamma(W\smallsetminus M_0,W\smallsetminus M;\pi)
		\end{equation*}
		is homotopic to the scanning map of McDuff \cite{MDCS} when $\pi:Y\rightarrow M$ is a trivial $X$-bundle. If we compose with a reflection $(t_1,\ldots,t_n;x)\mapsto (-t_1,\ldots,-t_{n-1},-t_n;x)$ in each fibre of $\dot\tau_\pi$ it is homotopic to the scanning map of B\"odigheimer and Madsen \cite{BMHQ} when $\pi$ is trivial and the extension due to Hesselholt \cite{HHD} for 
arbitrary $\pi$. 
	\end{remark}

        \begin{theorem}
                \label{thm:homotopy_equivalence}
                The scanning map $\gamma:E(M,M_0;\pi) \rightarrow \Gamma(W\smallsetminus M_0,W\smallsetminus M;\pi)$ is a weak homotopy equivalence if the pair $(M,M_0)$ or the typical fibre $X$ of $\pi$ is 0-connected. 
Moreover, if $X$ has the homotopy type of a CW-complex 
then it is a homotopy equivalence.
        \end{theorem}

        \begin{proof}
                By Remark \ref{rmk:scan_homotopic} we have 
a homotopy commutative diagram
                \begin{equation*}
                        \xymatrix{
                        E(M,M_0;\pi)\ar[r]^-{\gamma} &
                                \Gamma(W\smallsetminus M_0,W\smallsetminus M;\pi) \\
                        C(M,M_0;\pi)\ar[u]^{\simeq} \ar[ur]_{\simeq}
                        }
                \end{equation*}
                factoring the ordinary scanning map, 
which is a weak homotopy equivalence; see \cite {HHD} for the case 
with twisted coefficients. 
By proposition 2.11, the left hand map is a 
weak homotopy equivalence. 
Hence,   $\gamma$ is one too.
If $Y$ has the homotopy type of CW-complex so do 
all the three spaces in the diagram and hence the weak homotopy equivalences are homotopy equivalences.
	\end{proof}

		
	\subsection{Equivariance} 
	\label{sub:equivariance}
	

      Let $\Diff(M,\partial M)$ be the group of diffeomorphisms of $M$ which 
fix a collar of the boundary pointwise. 
There is an isomorphism $\Diff(M,\partial M) \rightarrow 
\Diff(W,W\smallsetminus M)$ extending diffeomorphisms by the identity 
on $W\smallsetminus M$. 
We will not distinguish between these groups and we understand the action 
of a diffeomorphism in $\Diff(M,\partial M)$ on $W$ to mean the action of its 
extension over $W\smallsetminus M$. 

Now we define actions on configuration spaces as follows. Let $\phi\in \Diff(M,\partial M)$, then $\phi\cdot(m_1,\ldots,m_k) \defeq (\phi(m_1),\ldots,\phi(m_k))$ for $(m_1,\ldots,m_k)\in\widetilde C_k(M)$ or $C_k(M)$ and $\phi\cdot f \defeq \phi\circ f\circ d\phi^{-1}$ for $f\in \widetilde E_k(M)$ or $E_k(M)$. 
With these actions the projections $p$ and $\widetilde p$ 
are $\Diff(M,\partial M)$-equivariant non-equivariant homotopy equivalences. 



This action of the diffeomorphism group  extends to 
configuration spaces with  labels when $Y = M\times X$ is a trivial bundle. 
However, in general diffeomorphisms of $M$ do not lift to automorphisms 
of the labelling bundle $Y$. 
Instead we need to consider the automorphism group of $Y$ 
and the induced diffeomorphisms on the base space $M$. 
	
	Let $\pi_1:Y_1\rightarrow M$ and $\pi_2:Y_2\rightarrow M$ 
be fibre bundles over $M$. 
We consider bundle isomorphisms $\alpha:Y_1\rightarrow Y_2$ that 
induce diffeomorphisms $\alpha_M$ on the base 
and denote isomorphic bundles by $Y_1\cong Y_2$. 
If the isomorphism induces the identity on the base 
we call it an equivalence and denote equivalent bundles by $Y_1\equiv Y_2$.
	
	\begin{lemma}
		Let $\Aut(\pi)$ denote the group of 
automorphisms of $\pi:Y\rightarrow M$ that restrict to diffeomorphisms of $M$.
Then the homomorphism $\Aut(\pi)\rightarrow \Diff(M)$, $\alpha\mapsto\alpha_M$ has image
		\begin{equation*}
			\Diff(M;\pi)\defeq\{\beta\in\Diff(M):(\beta_M)^*Y\equiv Y\}.
		\end{equation*}
	\end{lemma}
	
	\begin{proof}
		Let $\beta\in\Diff(M;\pi)$ and $\tilde\beta:\beta^*Y\rightarrow Y$ be the map completing the pullback square. Let $f:Y\rightarrow \beta^*Y$ be a bundle equivalence, then $\tilde\beta\circ f$ is an automorphism of $Y$ 
whose image under $\rho$ is $\beta$.
Conversely let $\alpha\in\Aut(\pi)$, and let $\widetilde{(\alpha_M)}:(\alpha_M)^*Y\rightarrow Y$ 
be the map completing the pullback square, then 
$\alpha^{-1}\circ\widetilde{(\alpha_M)}:(\alpha_M)^*Y\rightarrow Y$ 
is an equivalence.
	\end{proof}
	
	\begin{remark}
		In general, $\Diff(M;\pi)\subsetneq \Diff(M)$. For example, consider the Hopf bundle over $S^2$ and the antipodal map on the base. The Chern class of the Hopf bundle is $+1$ but the Chern class of the pullback is $-1$ so they cannot be equivalent. 
However, when $Y$ is trivial or some other natural bundle such 
as a tangent bundle, $\Diff(M;\pi) = \Diff(M)$ and there is a lift 
$\Diff(M)\rightarrow \Aut(\pi)$. 
For any bundle, $\Diff(M;\pi)$ contains $\Diff_0(M)$, 
the connected component of the identity, since homotopic maps induce equivalent pullbacks.
	\end{remark}
	

Let $\Aut^o(M,M_0\cup\partial M;\pi)\subset\Aut(\pi)$ denote the subgroup 
of bundle automorphisms of $Y$ which preserve the zero section and 
restrict to diffeomorphisms in $\Diff(M,M_0\cup\partial M)$. Here we assume as is standard
that the diffeomorphisms fix a collar in $M$ of $M_0 \cup \partial M$. 
	Let $\alpha\in\Aut^o(M,M_0\cup\partial M;\pi)$, then we define 
actions as follows:
	\begin{itemize}
		\item $\alpha\cdot(\mathbf{m};\mathbf{x})\defeq (\alpha_M(\mathbf{m});\alpha\circ\mathbf{x}\circ\alpha_M ^{-1})$ for $(\mathbf{m};\mathbf{x})\in C_k(M;\pi)$
		\item $\alpha\cdot(f;s)\defeq (\alpha_M\circ f\circ d\alpha_M^{-1}; \alpha \circ s\circ\alpha_M^{-1})$ for $(f;s)\in E_k(M;\pi)$
		\item $\alpha\cdot\sigma\defeq (d\alpha_M\wedge_W\alpha)\circ \sigma\circ \alpha_M ^{-1}$ for 
$\sigma\in\Gamma(W\smallsetminus M_0,W\smallsetminus M^+;\pi)$.
	\end{itemize}
	The action on an equivalence class in $C(M,M_0;\pi)$ or $E(M,M_0;\pi)$ is given by the equivalence class of the image of any representative under the action. This is 
well-defined as $\alpha$ preserves the zero section.
	
	
	\begin{proposition}
		The projection $p:E(M,M_0;\pi)\rightarrow C(M,M_0;\pi)$ is $\Aut^o(M,M_0\cup\partial M;\pi)$-equivariant.
	\end{proposition}
	
	\begin{proof}
		It follows immediately from the definitions that the projections $E_k(M;\pi)\rightarrow C_k(M;\pi)$ are equivariant. We obtain the result by noting that the actions are well behaved with respect to equivalence classes.
	\end{proof}
	
	\begin{remark}
		Although the projection is equivariant, the homotopy inverse  
$\sigma :C(M,M_0;\pi)\rightarrow E(M,M_0;\pi)$ defined by the exponential map and open 
$\varepsilon$-disks around each particle depends on choosing a metric and is not 
preserved under diffeomorphism.
	\end{remark}
	
	\begin{theorem}
		\label{thm:equivariance}
		The scanning map $\gamma:E(M,M_0;\pi) \rightarrow \Gamma(W\smallsetminus M_0,W\smallsetminus M;\pi)$ is $\Aut^o(M,M_0\cup\partial M;\pi)$-equivariant.
	\end{theorem}
	
	\begin{proof}
		Consider the scanning map $\gamma^+_k:
E_k(M;\pi)\rightarrow \Gamma(W\smallsetminus M_0,W\smallsetminus M^+;\pi)$ 
for some $k\ge 0$. Choose a point $z\in W\smallsetminus M_0$, 
an automorphism $\alpha\in\Aut^o(M,M_0\cup\partial M;\pi)$ 
and a tubular configuration $(f;s)$ with underlying 
configuration $(\mathbf{m};\mathbf{x})$. 
We want to show 
$$
(\alpha\cdot\gamma_k^+(f;s))(z) = \gamma_k^+(\alpha\cdot(f;s))(z).
$$
Suppose $\alpha_M^{-1}(z)\in\im(f_i)$ for some $i$. By the definition of the action on $\Gamma(W\smallsetminus M_0,W\smallsetminus M^+;\pi)$ we have
		\begin{equation*}
			(\alpha\cdot\gamma_k^+(f;s))(z)
			=\left(d_{\alpha_M^{-1}(z)}\alpha_M
			\circ d_{f_i^{-1}\circ \alpha_M^{-1}(z)}f_i
			\circ \phi_{f,\alpha_M^{-1}(z)}
			\circ f_i^{-1} \circ \alpha_M^{-1}(z),
			\alpha\circ s\circ\alpha_M^{-1}(z)\right).
		\end{equation*}
		We calculate $\gamma_k^+(\alpha\cdot(f;s))(z)$ componentwise: 
		\begin{enumerate}[(i)]
			\item $(\alpha\cdot f_i)^{-1} = d_{m_i}\alpha_M\circ f_i^{-1}\circ \alpha_M^{-1}$
			\item by Lemma \ref{lem:canonical_iso} there is a commutative diagram
		\begin{equation*}
			\xymatrix@C=40pt{
			T_{m_i}M \ar[r]^-{\phi_{f,\alpha_M^{-1}(z)}} \ar[d]_{d_{m_i}\alpha_M}
				& T_{f_i^{-1}\circ \alpha_M^{-1}(z)}T_{m_i}M \ar[d]^{d_{f_i^{-1}\circ\alpha_M^{-1}(z)}d_{m_i}\alpha_M}\\
				T_{\alpha_M(m_i)}M \ar[r]_-{\phi_{\alpha\cdot f,z}} & T_{(\alpha\cdot f_i)^{-1}(z)}T_{\alpha_M(m_i)}M
			}
		\end{equation*}
		so we have $\phi_{\alpha\cdot f,z}
		=d_{f_i^{-1}\circ\alpha_M^{-1}(z)}d_{m_i}\alpha_M
		\circ\phi_{f,\alpha_M^{-1}(z)} \circ\left(d_{m_i}\alpha_M\right)^{-1}$
		\item $d_{(\alpha\cdot f_i)^{-1}(z)}(\alpha\cdot f_i) = d_{\alpha_M^{-1}(z)}\alpha_M\circ d_{f_i^{-1}\circ\alpha_M^{-1}(z)}f_i \circ\left(d_{f_i^{-1}\circ\alpha_M^{-1}(z)}d_{m_i}\alpha_M\right)^{-1}$.
		\end{enumerate}
		Putting these together we have
		\begin{equation*}
			\gamma_k(\alpha\cdot(f;s))(z)
			=\left(d_{\alpha_M^{-1}(z)}\alpha_M
			\circ d_{f_i^{-1}\circ \alpha_M^{-1}(z)}f_i
			\circ \phi_{f,\alpha_M^{-1}(z)}
			\circ f_i^{-1} \circ \alpha_M^{-1}(z),
			\alpha\circ s\circ\alpha_M^{-1}(z)\right).
		\end{equation*}
In our computations above we have implicitly assumed that none of the factors of
$(f;s)$ have midpoint in the collar  of $\partial M_0$ in $M$.
We leave it to the reader to check that the multiplicative factor $\exp (\frac 1 {
\delta (m_i)})$ does not interfere with the equivariance argument as 
our diffeomorphisms fix the collar by assumption.

The homotopy 
$\Gamma (W\smallsetminus M_0, W\smallsetminus M^+; \pi ) \to 
\Gamma (W\smallsetminus M_0, W\smallsetminus M; \pi )$ is invariant 
under $\alpha$, again
since automorphisms in $\Aut^o(M,M_0\cup\partial M;\pi)$ fix the chosen collar. 
To complete the proof recall that the scanning map $\gamma$ is defined by 
$\rho \circ \gamma_k^+$ for any representative of a class 
and the actions respect the equivalence relation.
	\end{proof}
	
       Combining Theorems \ref{thm:equivariance} and \ref{thm:homotopy_equivalence} gives the following corollary.

        \begin{corollary}
		Let $(M, M_0)$ or a typical fibre $X$ of $\pi$ be 0-connected.
                If a topological group $G$ acts on $M$ and $\pi$ via 
a homomorphism $G\rightarrow \Aut^o(M,M_0\cup\partial M;\pi)$, then there are weak homotopy equivalences of the reduced Borel constructions
                \begin{equation*}
                        EG\ltimes_GC(M,M_0;\pi) \simeq EG\ltimes_GE(M,M_0;\pi) \simeq EG\ltimes_G\Gamma(W\smallsetminus M_0,W\smallsetminus M;\pi)
                \end{equation*}
                which are homotopy equivalences if $Y$ has the homotopy type of a CW-complex.
        \end{corollary}

        \begin{remark}
                Let $Y = M \times X$ be a trivial bundle and let $G$ 
be a compact Lie group acting on $M$ through isometries  and  on $X$ 
fixing the base-point. Then the ordinary scanning map
$$
\gamma \circ \sigma : C(M, M_0; X) \longrightarrow 
\Gamma (W \smallsetminus M_0, W\smallsetminus M; X)
$$
is $G$-equivariant. 
        \end{remark}

	
	
	\section{Stable splittings} 
	\label{sec:stable_splittings}
	
	In this section we construct Snaith type  splittings 
via  diffeomorphism
equivariant maps.  As a result we generalise and strengthen the results 
of B\"odigheimer and Madsen \cite{BMHQ}. As an application we 
also derive an equivariant 
stable split 
injection for configuration spaces and diffeomorphisms of manifolds 
with marked points. This generalises some 
results in \cite{BT}.

{\it Constructions and assumptions:} 
Let $M_0 \subset M \subset M^+ \subset W$ and $\pi:Y \to X$ be as before.
For each $k\ge 1$ let $C^k = C^k(M,M_0;\pi)$ be the filtrations of $C=C(M,M_0;\pi)$.
For each $k\ge 1$ let $D^k= D^k (M, M_0; \pi) \defeq C^k/C^{k-1}$ be the filtration quotient and let $V\defeq \bigvee_{i\ge 1}D^i$ be the wedge sum with $k$th filtration $V^k\defeq \bigvee_{i=1}^kD^i$.

Throughout subsections 4.1-4.4 we will assume that 
$(M,M_0)$ or the typical fibre $X$  of $\pi$ is $0$-connected.
Then also $C, C^k,  D^k, V, V^k $ are all path connected as well.
We also note that if $X$ has the homotopy type of a CW-complex then 
so do all these spaces.

Throughout this section we let $G$ be a topological group which 
acts on $M$ and $\pi$ via a homomorphism $G\rightarrow 
\Aut^o(M,M_0\cup\partial M;\pi)$. 
We will suppress this homomorphism in our notation.
Note that the $G$ action on $C$ induces an action on $C^k$, $D^k$ and hence on $V^k$.

		
	
	\subsection{Equivariant splittings in homology} 
	\label{sub:equivariant_splittings_in_homology}
	
	Before we construct equivariant stable splittings for mapping spaces we prove the following weaker statement in homology in virtue of the  simplicity 
of its proof.
	
	\begin{theorem}
		There exist  a $\pi_0G$-equivariant isomorphism
		\begin{equation*}
			\widetilde H_*(C(M,M_0;\pi))\xlongrightarrow{\cong}  \bigoplus_{k\ge 1} \widetilde H_*(D^k).
		\end{equation*}
	\end{theorem}
	
	\begin{proof}
		Given a pointed space $A$ we denote the infinite 
symmetric product of $A$ by $\SP(A)$. 
Elements of $\SP(A)$ are finite formal sums $\sum k_ia_i$ 
where $a_i\in A$ and $k_i\in \mathbb{N}$ for each $i$ with the one relation 
that makes the base point the zero of the monoid, see \cite{DTQ}. 
		Let $\xi_\alpha$ be  a configuration in $C$ 
and write it as a formal sum 
$\xi_\alpha = \sum_{i\in\alpha}m_ix_i$ where $\alpha$ is a finite set. Given a subset $\beta\subseteq\alpha$ of size $k$ there is a subconfiguration $\xi_\beta\defeq \sum_{i\in\beta}m_ix_i\in C^k$. Let $\bar\xi_\beta$ denote its image under the composition $C^k\rightarrow D^k\hookrightarrow V$. Define a power set map
		\begin{align*}
			\sigma:C &\longrightarrow\SP(V)\\
			\xi_\alpha &\longmapsto \sum_{\beta\subseteq\alpha}\xi_\beta
		\end{align*}
		and extend it to a map $\sigma:\SP(C)\rightarrow\SP(V)$ by $\sum k_\alpha\xi_\alpha\mapsto \sum k_\alpha\sigma(\xi_\alpha)$. 
This is a pointed $G$-map that respects the filtrations and restricts to a pointed $G$-map $\SP(C^k)\rightarrow\SP(V^k)$.
		
		The infinite symmetric product functor sends cofibrations to quasifibrations so for each $k\ge 1$ we have a commutative diagram
		\begin{equation*}
			\xymatrix{
			\SP(C^{k-1}) \ar[d]\ar[r]^\sigma & \SP(V^{k-1}) \ar[d] \\
			\SP(C^k) \ar[d]\ar[r]^\sigma & \SP(V^k) \ar[d] \\
			\SP(C^k\big/C^{k-1}) \ar[r]^{\sigma=\id} & \SP(V^k\big/V^{k-1})
			}
		\end{equation*}
		in which the vertical sequences are quasifibrations and the lower horizontal map is the identity $\SP(D^k)\rightarrow \SP(D^k)$. Starting with $C^1 = V^1$ we proceed by induction on $k$ to see that $\sigma$ induces isomorphisms
		\begin{equation*}
			\pi_*\SP(C(M,M_0;\pi))\xlongrightarrow{\cong}\pi_*\SP(\bigvee_{k\ge 1}D^k).
		\end{equation*}
		By the Dold-Thom theorem \cite{DTQ} this is precisely the homology isomorphism we require. To finish the proof note that the connected component of the identity $G_0\subset G$ acts trivially on homology so the isomorphism is equivariant under the actions of $\pi_0G = G/G_0$.
	\end{proof}
	
Combining this with the results in section 3 we have proved the 
following corollary.

\begin{corollary}
	There exists a $\pi_0 G$-equivariant isomorphism
$$
\widetilde H_* ( \Gamma (W \smallsetminus M_0, W \smallsetminus M; \pi)) 
\xlongrightarrow{\cong}  \bigoplus_{k\ge 1} \widetilde H_*(D^k).
\qed
$$
\end{corollary}

	The above  proof serves as a model for the proofs of the stronger
theorems 
in the following sections. 
The main difficulty in each case is to construct a suitable 
$G$-equivariant power set map. 
	
	
	\subsection{Stable splittings of reduced Borel constructions} 
	\label{sub:stable_splittings_of_reduced_borel_constructions}
	
	B\"odigheimer and Madsen's power set map for Borel constructions
in \cite{BMHQ} relies crucially on 
the fact that the groups they consider are compact of Lie type.
In particular they use that one can embed the configuration 
space $C(M)$ $G$-equivariantly into a finite dimensional 
$G$-representation. 
Instead we use here suitably chosen models for $EG$ 
which allow us to define the power map on Borel constructions for any 
$G$ acting on $M$  smoothly and on $\pi$. 
This provides a geometric model of the power set map.
	
	\begin{theorem}
		There is a weak homotopy equivalence of suspension spectra
		\begin{equation*}
			\Sigma^\infty\left(EG\ltimes_G C(M,M_0;\pi)\right)\longrightarrow \Sigma^\infty\bigvee_{k\ge1} EG\ltimes_G D^k.
		\end{equation*}
	\end{theorem}
	
	\begin{proof}
		As a consequence of Whitney's embedding theorem,
the embedding space $\Emb(M;\mathbb{R}^\infty)$ is contractible. 
It admits a free $\Diff(M)$-action. 
Choose a model for $EG$ and replace it with 
$EG\times \Emb(M;\mathbb{R}^\infty)$ equipped with the diagonal 
action where $g\cdot h \defeq h\circ (g^{-1})$ for $g\in G$ and an embedding $h:M\rightarrow\mathbb{R}^\infty$. This is again a model for $EG$ and we denote elements of this $EG$ by pairs $(u,h)$.
		
		Given an embedding $h:M\rightarrow\mathbb{R}^\infty$ 
we define an associated 
embedding $\widehat h:\coprod C_k(M)\rightarrow \mathbb{R}^\infty$
as follows. 
For each $k \geq 1$ and any $N$, $C_k (\mathbb R^N)$ is a finite dimensional,
smooth manifold.   Thus  there is
an embedding $C_k (\mathbb R^N) \hookrightarrow
\mathbb R ^{L_N} $ for some $L_N$. Considering $C_k (\mathbb R^N) $ 
as a submanifold of $C_{k} (\mathbb R^{N+1}) $ we can extend 
this embedding to an embedding of $C_{k} (\mathbb R^{N+1}) \hookrightarrow
\mathbb R^{L_{N+1}}$ for some $L_{N+1} \geq L_N$. Proceeding like this, 
in the limit we will have constructed an embedding 
$$
\mu _k: C_k( \mathbb R^\infty) \hookrightarrow \mathbb R^\infty.
$$
Choosing an injection 
$\coprod _{k\geq 1} \mathbb R^\infty \hookrightarrow \mathbb R^\infty$
we can construct from  $\mu_k$ an injection 
$$
\mu: 
\coprod _{k\geq 1} C_k (\mathbb R^\infty) \hookrightarrow \mathbb R^\infty.
$$
Given now any embedding $h: M \hookrightarrow \mathbb R^\infty$, 
the embedding $\widehat h$ 
is defined as the composition
		\begin{equation*}
			\widehat h: \coprod_kC_k(M) \xlongrightarrow{C_kh} \coprod_k C_k(\mathbb{R}^\infty) \xlongrightarrow{\mu} \mathbb{R}^\infty.
		\end{equation*}
The action of $G$ on $M$ cannot be extended to an action on 
$\mathbb R^\infty$, and in particular 
$\widehat h$ cannot be made  $G$-equivariant. Nevertheless, the  
construction is 
$G$-equivariant in the sense that  
$$
g. \widehat
h = \widehat h \circ {g^{-1}} = \widehat {h \circ g^{-1}} = \widehat {g.h}. 
$$
We can now define the power-set map in this setting.
Choose a configuration $\xi_\alpha = \sum_{i\in\alpha}m_ix_i\in C$ 
and define subconfigurations $\xi_\beta\in C^k$ and $\bar\xi_\beta\in V$ 
as in the previous section. 
Let $m_\beta\defeq\sum_{i\in\beta}m_i\in C_k(M)$ 
be the associated unlabelled subconfiguration and define
		\begin{align*}
			\sigma:EG\ltimes_G C(M,M_0;\pi) &\longrightarrow C(\mathbb{R}^\infty;EG\ltimes_G V) \\
			(u,h,\xi_\alpha) &\longmapsto \sum_{\beta\subseteq\alpha}\widehat h(m_\beta)(u,h,\bar\xi_\beta)
		\end{align*}
		where the target is the configuration space of 
unordered particles in $\mathbb{R}^\infty$ with labels in the 
trivial $EG\ltimes_G V$-bundle. 
The formula for  $\sigma$  gives a well defined map 
$EG\times C\rightarrow C(\mathbb{R}^\infty;EG\times V)$. 
To see that it is well defined on the half smash product note that the 
base-point in $C$ is the empty configuration $\xi_\emptyset$ 
which has no subconfigurations so $\sigma$ maps $(u,h,\xi_\emptyset)$ 
to the empty configuration in the target for any $(u,h)\in EG$.
To check that $\sigma$ is $G$-equivariant
recall that $g \cdot\xi_\alpha =\sum g(m_i)\left(g\circ x_i\circ g^{-1}\right)$.
The $G$-action on the target is trivial on $\mathbb{R}^\infty$ 
and given by the diagonal action on the labels.
Hence,
	\begin{align*}
		\sigma (g\cdot(u,h,\xi_\alpha) ) 
&= \sum_{\beta\subseteq\alpha} \widehat{h\circ g^{-1}} 
(\sum_{i\in\beta} g(m_i) )
\left(g\cdot u, 
h\circ g^{-1},\overline{g \cdot\xi}_\beta\right)\\
		g\cdot\sigma(u,h,\xi_\alpha) &= \sum_{\beta\subseteq\alpha}\widehat h(m_\beta)(g\cdot u, h\circ g^ {-1}, g \cdot\bar\xi_\beta )
	\end{align*}
	By definition $g \cdot\bar\xi_\beta = \overline{g\cdot\xi}_\beta$ and 
$$
\widehat{h\circ g^{-1}} (\sum_{i\in\beta} g(m_i)) 
= \mu (\sum_{i\in\beta}h\circ g^{-1}\circ g(m_i)) 
=\mu (\sum _{i \in \beta} h(m_i)) 
= \widehat h(m_\beta).
$$ 
So $\sigma$ is well defined.
	
	Composing  the power set map with the scanning map yields the map
	\begin{equation*}
		EG\ltimes_G C\xlongrightarrow{\sigma} C(\mathbb{R}^\infty;EG\ltimes_G V) \longrightarrow \Omega^\infty\Sigma^\infty (EG\ltimes_G V)
	\end{equation*}
	which respects the filtrations of $C$ and $V$. The cofibration sequences of the filtrations of $C$ and $V$ induce fibration sequences of suspension spectra which fit into commutative diagrams
	\begin{equation*}
		\xymatrix{
		\Sigma^\infty\left(EG\ltimes_G C^{k-1}\right) \ar[d]\ar[r]^p
			& \Sigma^\infty\left(EG\ltimes V^{k-1}\right) \ar[d]\\
		\Sigma^\infty\left(EG\ltimes_G C^k\right) \ar[d]\ar[r]^p
			& \Sigma^\infty\left(EG\ltimes V^k\right) \ar[d]\\
		\Sigma^\infty\left(EG\ltimes_G C^k\big/C^{k-1}\right) \ar[r]^{p \simeq \id}
			& \Sigma^\infty\left(EG\ltimes V^k\big/V^{k-1}\right)
		}
	\end{equation*}
	for each $k\ge 1$ and where $p$ is the adjoint of the composition above. The lower horizontal map is homotopic to the identity since it is the adjoint of the inclusion $EG\ltimes_G D^k\hookrightarrow \Omega^\infty\Sigma^\infty\left(EG\ltimes_G D^k\right)$. Starting with $C^1 = V^1$ we proceed by induction on $k$ to obtain a weak homotopy equivalence of suspension spectra and the result follows from the observation that $EG\ltimes_G \bigvee_{k\ge1}D^k = \bigvee_{k\ge1}EG\ltimes_G D^k$.
	\end{proof}

	
	\subsection{Equivariant stable splittings} 
	\label{sub:equivariant_stable_splittings}

	The purpose of this section is to construct a $G$-equivariant 
power-set map which  induces stable $G$-equivariant splittings  for 
configuration spaces. 	
It will be completely natural and avoid having to choose an embedding $\mu$ 
as in the proof of Theorem 4.3. It also gives a stronger result.

	While the stable splitting in the previous section made use of the 
configuration space $C(\mathbb{R}^\infty;A)$ as a 
model for the free infinite loop space $Q(A)\defeq\rlim\,\Omega^n\Sigma^nA$, 
here we will make use of 
the $\Gamma^+$-construction of Barratt and Eccles.  We will also need to
replace 
the $\Sigma_k$-orbits of $\widetilde C_k(M; \pi)$ by the homotopy  
$\Sigma _k$-orbits. 
Recall the following constructions from \cite {BEG}.

	For a discrete group $H$, let $W.H$ denote its homogeneous 
bar construction  and let $EH\defeq|W.H|$ be its 
geometric realisation. 
Group homomorphisms $H\rightarrow K$ induce continuous maps $EH\rightarrow EK$,
and in particular the inclusion $\Sigma_{k-1}\hookrightarrow \Sigma_k$ 
induces an inclusion map $E\Sigma_{k-1}\hookrightarrow E\Sigma_{k}$ 
with a right inverse induced by the reduction map 
$R:\Sigma_k\rightarrow \Sigma_{k-1}$ as defined in \cite {BEG}.
	Given a well-pointed space $A$, the $\Gamma^+$ construction on $A$ is defined as
	\begin{equation*}
		\Gamma^+(A)\defeq \left(\coprod_{k\ge 0}E\Sigma_k\times_{\Sigma_k} A^k\middle)\right/\sim
	\end{equation*}
	where $(w;a_1,\ldots,a_k)\sim (R(w);a_1,\ldots,a_{k-1})$ if $a_k = *$.
	$\Gamma^+$ is an endofunctor on the category of well-pointed 
topological spaces. For any space $A$ there is an 
inclusion map $\iota_A:A\hookrightarrow\Gamma^+(A)$ 
identifying $A$ with $\Sigma_1\times_{\Sigma_1}A\subset\Gamma^+(A)$ 
and a structure map $h_A:\Gamma^+\Gamma^+(A)\rightarrow \Gamma^+(A)$ 
making the triple $(\Gamma^+,\iota,h)$ into a monad. 
The spaces $\Gamma^+(A)$ are the free $\Gamma^+$-algebras.
	
	We  define the 
Borel configuration space of $k$ unordered 
particles in $M$ with twisted labels in $\pi$ as
	\begin{equation*}
		\mathscr{C}_k(M;\pi)\defeq E\Sigma_k\times_{\Sigma_k}\widetilde C_k(M;\pi).
	\end{equation*}
	The Borel configuration space with particles vanishing in $M_0$ is then
	\begin{equation*}
		\mathscr{C}(M,M_0;\pi)\defeq \left(\coprod_{k\ge 0} \mathscr{C}_k(M;\pi)\middle)\right/\sim
	\end{equation*}
	where $(w;m_1,\ldots,m_k;x_1,\ldots,x_k)\sim (R(w);m_1,\ldots,m_{k-1};x_1,\ldots,x_{k-1})$ if $m_k\in M_0$ or $x_k=o|_{m_k}$.
	The Borel configuration spaces admit filtrations $\mathscr{C}^k\defeq \coprod_{i=0}^k\mathscr{C}_k\big/\sim$ analogous to the ordinary configuration spaces. We make use of the filtrations quotients $\mathscr{D}^k\defeq \mathscr{C}^k\big/\mathscr{C}^{k-1}$. 
Consider the sequence of spaces $\widetilde D^k = \widetilde C^k(M;\pi)\big/\equiv$ where $(m_1,\ldots,m_k;x_1,\ldots,x_k)\equiv *$ 
if $m_i\in M_0$ or $x_i=o|_{m_i}$ for some $i$. 
These spaces are  the filtration quotients of the ordered 
configuration space with vanishing particles. 
The filtration quotients of the Borel configuration spaces can 
then be identified  as 
$\mathscr{D}^k=E\Sigma_k\ltimes_{\Sigma_k}\widetilde D^k$.	
	
	\begin{proposition}
		\label{prop:Borel_projection_equivalence}
The projection $P: \mathscr{C}(M,M_0;\pi)\rightarrow C(M,M_0;\pi)$ 
and the restriction to the filtration quotients 
$P: \mathscr{D}^k\rightarrow D^k$ are  weak homotopy equivalences.
		\end{proposition}
	
	\begin{proof}
		The projection $P:E\Sigma_k\times\widetilde C_k(M;\pi)
\rightarrow\widetilde C_k(M;\pi)$ is a $\Sigma_k$-equivariant 
homotopy equivalence. As $\Sigma _k$ acts freely on the source and target, 
$P$ induces 
a homotopy equivalence on orbit spaces 
$\mathscr{C}_k(M;\pi)\xrightarrow{\simeq} C_k(M;\pi)$. 
Together they give a well defined map $P:
\mathscr{C}(M,M_0;\pi)\rightarrow C(M,M_0;\pi)$ of filtered spaces.
		
		Let $\widetilde B_{C_k}(M,M_0;\pi)\subset\widetilde C_k(M;\pi)$
 be the subspace of ordered configurations such that for 
some $i$ either $m_i\in M_0$ or $x_i=o|_{m_i}$.  Recall from the proof of Proposition \ref{prop:proj_weak_equiv} that the filtration $C^k(M,M_0;\pi)$ is obtained as the pushout of the diagram $C^{k-1}(M,M_0;\pi)\leftarrow B_{C_k}(M,M_0;\pi)\rightarrow C_k(M;\pi)$ where $B_{C_k}(M,M_0;\pi) = \widetilde B_{C_k}(M,M_0;\pi)\big/\Sigma_k$. Similarly the filtrations of the Borel configuration spaces are obtained from the previous filtration as the pushout along a subspace $B_{\mathscr{C}_k}(M,M_0;\pi)\defeq E\Sigma_k\times_{\Sigma_k}\widetilde B_{C_k}(M,M_0;\pi)$. The projection $\pi$ induces a homotopy equivalence $B_{\mathscr{C}_k}(M,M_0;\pi)\xrightarrow{\simeq} B_{C_k}(M,M_0;\pi)$ and a map of homotopy pushout squares. The result then follows by induction as in the proof of Proposition \ref{prop:proj_weak_equiv}. The result for the filtration quotients is proved similarly.
	\end{proof}
	
	\begin{remark}
		The maps in the previous proposition are $G$-equivariant.
	\end{remark}
	
	\begin{proposition}
		There is a $G$-equivariant power-set map $\sigma: 
\mathscr{C}(M,M_0;\pi)\rightarrow \Gamma^+(V)$ of filtered spaces.
	\end{proposition}
	
	\begin{proof}
		Fix $k\ge 0$. Let $\alpha$ be an ordered set of cardinality $k$,
for each $0\leq i\leq k$ let $\mathscr{P}_i(\alpha)$ be the set of (unordered)
subsets of $\alpha$ of size $i$. As $\alpha$ is ordered 
$\mathscr P _i (\alpha)$ has an induced lexicographical ordering and 
we may list its elements  as:
 $\beta_1,\ldots,\beta_{k\choose i}$. 

Let $\xi_\alpha= \sum_{j\in\alpha} m_jx_j$ 
be an ordered configuration in $\widetilde C_k(M;\pi)$. 
For each $\beta\in \mathscr{P}_i(\alpha)$ let $\xi_\beta=\sum_{i\in\beta}m_jx_j\in C_i(M;\pi)$ be an unordered subconfiguration of $\xi_\alpha$ and let $\bar\xi_\beta$ be its image under the composition $C_k\rightarrow C^k\rightarrow D^k\rightarrow V$. For each $i$ define a map
		\begin{align*}
			\widetilde C_k(M;\pi) &\longrightarrow V^{k\choose i}\\
			\xi_\alpha &\longmapsto \left(\bar\xi_{\beta_1},\ldots,\bar\xi_{\beta_{k\choose i}}\right).
		\end{align*}
		
$\Sigma_k$ naturally acts on $\mathscr P _i(\alpha)$ since permuting the 
elements in $\alpha$  maps a subset  of size $i$ to another 
subset of size $i$. 
This action defines a homomorphism $\Sigma_k\rightarrow \Sigma_{k\choose i}$.
With this action, the above map is $\Sigma_k$-equivariant. 
For each $k$ the direct sum of these homomorphisms defines a homomorphism $\phi_k:\Sigma_k\rightarrow \Sigma_{k\choose 0}\oplus\Sigma_{k\choose 1}\oplus\cdots\oplus\Sigma_{k\choose k}\subset\Sigma_{2^k}$. Together with the maps $\widetilde C_k(M;\pi)\rightarrow V^{k\choose i}$ this defines
		\begin{equation*}
			\sigma_k:\mathscr{C}_k(M;\pi) = E\Sigma _k \times _{\Sigma_k}
\widetilde C_k (M; \pi) 
 \longrightarrow E\Sigma_{2^k}\times_{\Sigma_{2^k}}V^{2^k}.
		\end{equation*}
		
We will now show that the $\sigma_k$'s fit together to define a map 
$\sigma:  \mathscr C(M;\pi) \to \Gamma ^+ (V)$. 
To see that $\sigma$ respects the equivalence relations on each side, 
choose two equivalent Borel configurations. 
We may assume they are of the form $(w;m_1,\ldots,m_k;x_1,\ldots,x_k)$ and its
 reduction $(R(w);m_1,\ldots,m_{k-1};x_1,\ldots,x_{k-1})$. 
Since $x_k=o|_{m_k}$ or $m_k\in M_0$, any subconfiguration 
containing $(m_k;x_k)$ is mapped to the base-point in the appropriate 
copy of $V\subset V^{2^k}$. The defining relation of 
$\Gamma^+$ identifies a point in $V^n$ with one in $V^{n-1}$ 
if a coordinate lies in the base-point of $V$. 
Precisely $2^{k-1}$ subconfigurations contain the particle $(m_k;x_k)$.
 So noting that the square
	\begin{equation*}
		\xymatrix{
		\Sigma_k\ar[r]^{\phi_k}\ar[d]_R & \Sigma_{2^k} \ar[d]^{R^{2^{k-1}}}\\
		\Sigma_{k-1} \ar[r]_{\phi_{k-1}} & \Sigma_{2^{k-1}}
		}
	\end{equation*}
	commutes up to an isomorphism $\Sigma_{2k}\rightarrow \Sigma_{2k}$ and iterating the relation $2^{k-1}$ times we have
	\begin{equation*}
		\sigma_k(w;m_1,\ldots,m_k;x_1,\ldots,x_k)\sim 
\sigma_{k-1}(R(w);m_1,\ldots,m_{k-1};x_1,\ldots,x_{k-1}).
	\end{equation*}
	
	To see that $\sigma$ is a map of filtered spaces, note that 
by definition $\sigma$ maps 
 $\mathscr{C}^k$ to $\Gamma^+(V_k)\subset\Gamma^+(V)$
where  $V_k\defeq \bigvee_{i=1}^kD^k$.
Finally, the construction is natural with respect to any 
inclusions of the underlying fibration $\pi$ into another fibration,  
and in particular it is  $G$-equivariant.
	\end{proof}
	
	\begin{theorem}
There is a stable  homotopy equivalence
		\begin{equation*}
			\mathscr{C}(M,M_0;\pi)\simeqs\bigvee_{k\ge 1}D^k(M,M_0;\pi)
		\end{equation*}
		induced by a $G$-equivariant map $\bar \sigma:\Gamma^+(\mathscr{C})\rightarrow \Gamma^+(V)$.
	\end{theorem}
	
	\begin{proof}
		Let $\bar \sigma$ denote the composition $\Gamma^+(\mathscr{C}(M,M_0;\pi)) \xlongrightarrow{\Gamma^+\sigma} \Gamma^+\Gamma^+(V) \xlongrightarrow{h_V} \Gamma^+(V)$ 
and note that it restricts to a map of the filtrations and 
filtration quotients. 
$\bar \sigma$ is $G$-equivariant as both the $\Gamma ^+$-construction on a map 
and the construction $h_V$ are natural.
		We note now that the image of 
$\mathscr{D}^k$ under $\sigma$ is contained in 
$E\Sigma_1\times_{\Sigma_1}D^k = D^k$ and hence $\sigma$ restricted to $
\mathscr D ^k$ factors as 
$\mathscr{D}^k\rightarrow D^k\xhookrightarrow{\iota_{D^k}}\Gamma^+(D^k)$ 
where the first map is the weak homotopy equivalence in Proposition \ref{prop:Borel_projection_equivalence}.
Thus on the filtration quotients $\bar \sigma$ factors as the composition
		\begin{equation*}
			\Gamma^+(\mathscr{D}^k) \xlongrightarrow{\Gamma^+ P} 
\Gamma^+(D^k) 
\xlongrightarrow{\Gamma^+\iota_{D^k}} \Gamma^+\Gamma^+(D^k) 
\xlongrightarrow{h_{D^k}} \Gamma^+(D^k).
		\end{equation*}
		The first map is a homotopy equivalence 
since $\Gamma^+$ is a homotopy functor. Also, 
the composition $h_{D^k}\circ\Gamma^+\iota_{D^k}$ is the identity 
by virtue of $\Gamma ^+$ being a monad.
It follows that $\bar \sigma:\Gamma^+(\mathscr{D}^k)\rightarrow \Gamma^+(D^k)$
is a homotopy equivalence.		

		The $\Gamma^+$ construction converts cofibrations to fibrations so for each $k$ we can form a commutative diagram in which the vertical sequences are fibrations.
		\begin{equation*}
			\xymatrix{
			\Gamma^+(\mathscr{C}^{k-1}) \ar[r]^{\bar \sigma}\ar[d] & \Gamma^+(V^{k-1}) \ar[d]\\
			\Gamma^+(\mathscr{C}^{k}) \ar[r]^{\bar \sigma}\ar[d] & \Gamma^+(V^{k}) \ar[d]\\
			\Gamma^+(\mathscr{D}^k) \ar[r]^{\bar \sigma} & \Gamma^+(D^k)
			}
		\end{equation*}
		By the above discussion the bottom arrow is always a homotopy equivalence. Since $\mathscr{C}^1=\mathscr{D}^1$ and $C^1=D^1=V^1$ we proceed by induction on $k$ to obtain a homotopy equivalence $\Gamma^+(\mathscr{C})\simeq \Gamma^+(V)$.

		By the results of  \cite{BEG2}, a  connected 
(or group-like) $\Gamma^+$-algebra is an infinite loop space 
and a morphism of such $\Gamma^+$-algebras is a map of infinite loop spaces. 
Moreover, the free $\Gamma^+$-algebras correspond to free infinite loop spaces if the underlying space is path connected and has the homotopy type of a CW-complex. By our assumptions on $(M,M_0)$ and $X$ 
all the spaces considered are path connected 
and have the homotopy types of CW-complexes. Furthermore, 
by definition $\bar \sigma$ and 
its restrictions are morphisms of $\Gamma^+$-algebras. 
Thus the statement of the theorem follows.
	\end{proof}

	
\subsection{Induced splittings for mapping spaces}

Combining the  stable splitting of Theorem 4.7 with the fact that the 
scanning map is equivariant, Theorem 3.11,  we immediately deduce the following result.

\vskip .1in
	\begin{theorem}	
 There is a zigzag of $G$-equivariant maps
		\begin{equation*}
			\Gamma(W\smallsetminus M_0,W\smallsetminus M;\pi)\longleftarrow E(M,M_0;\pi)\longrightarrow C(M,M_0;\pi) \longleftarrow \mathscr{C}(M,M_0;\pi) \longrightarrow \bigvee_{k\ge 1} D^k(M,M_0;\pi)
		\end{equation*}
		each of which is a stable homotopy equivalence. The zigzag induces a weak homotopy equivalence
		\begin{equation*}
			\Omega^\infty\Sigma^\infty \left(EG\ltimes_G \Gamma(W\smallsetminus M_0,W\smallsetminus M;\pi)\right) \simeq \prod_{k\ge 1} \Omega^\infty\Sigma^\infty \left(EG\ltimes_G D^k(M,M_0;\pi)\right).
		\end{equation*}
\qed
	\end{theorem}

We  now shift emphasis from configuration spaces to mapping spaces and adopt
the notation 
$\bar D^k$ for the $k$th filtration quotient of 
$$
C(M\smallsetminus M_0, \partial M \smallsetminus M_0; \pi).
$$
When $M$ is parallelisable and $\pi$ is a trivial bundle, the section space
in the above theorem is simply a mapping space.
But even when $M$ is not parallelisable, 
we can use configuration spaces with 
twisted coefficients in the normal bundle of
$M$ to \lq untwist' the section spaces so that they become again
mapping spaces.

	\begin{corollary}
For $N$ large enough, there is a  homotopy equivalence
		\begin{equation*}
			\Omega^\infty\Sigma^\infty\left(EG\ltimes_G\Map(M,M_0;\Sigma^NX,*)\right) \simeq \prod_{k\ge 1}\Omega^\infty\Sigma^\infty\left(EG\ltimes_G \bar D^k\right).
		\end{equation*}
	\end{corollary}
	
	\begin{proof}
Let $\nu$ be a bundle such that $TM\oplus \nu$ is a trivial bundle. (To 
construct such a $\nu$ note that
$M$ is compact and can be  embedded in $\mathbb{R}^N$ 
for some $N$ large enough.  $\nu$ may be taken to be the normal bundle 
of the embedding.)
Let $\pi:Y=\dot\nu\wedge_M(M\times X)\rightarrow M$ be the fibre bundle 
obtained by taking the fibrewise smash product of $X$ with the fibrewise 
one-point compactification of $\nu$, equipped with the canonical 
section at infinity.
Note that 
$$
\tau_\pi|_M = \dot TM\wedge_M\dot\nu\wedge_M(M\times X) \cong 
\Sigma^NX\times M
$$ 
The result is now a special case of the previous theorem.
	\end{proof}

\begin{example}
		If $M$ is stably parallelisable,
then for some $N\in\mathbb{N}$ there exists 
a trivial vector bundle $\nu$ such that $TM \oplus \nu$ is the trivial  $N$-plane bundle on $M$.
As $\nu$ is trivial, it admits 
an action of the full diffeomorphism group of $M$.
Then there is a homotopy equivalence
		\begin{equation*}
			\Omega^\infty\Sigma^\infty\left(E\Diff(M) \ltimes_{\Diff(M)}\Map(M;\Sigma^NX)\right) \simeq \prod_{k\ge 1}\Omega^\infty\Sigma^\infty\left(E\Diff(M) \ltimes_{\Diff(M)} \bar D^k\right).
		\end{equation*}

When $M$ is parallelisable  we may take $N$ equal to the dimension of $M$.
For  a compact Lie group $G$ and assuming further that $M$ 
is $G$-parallelisable, this recovers the main 
theorem of \cite{BMHQ}.
\end{example}

\begin{example}
When $M$ is the disc in $\mathbb R^n$ and $M_0$ empty,
 we see that the Snaith stable splitting of $\Omega ^n \Sigma^n (X)$ 
is equivariant under the action of $\Diff (\mathbb R^n)$.
\end{example}

\begin{example}
When $M$ is the $n$-sphere $S^n \subset \mathbb R^{n+1}$ 
and $M_0$ empty, then the normal bundle $\nu$ 
is the one dimensional trivial bundle. Thus for 
$\pi:Y=\dot\nu\wedge_M(M\times X) \simeq M\times \Sigma X \rightarrow M$  
$$
C(S^n, \emptyset ; \pi) \simeq \Map (S^n; \Sigma ^{n+1} X).
$$	
This provides a model for the free higher  loop space and
a stable splitting, both of which are $\Diff (S^n)$-equivariant. 
\end{example} 


\subsection {Stable injectivity for configurations and diffeomorphisms}

Let $M$ be as before, $M_0 = \emptyset$ 
and  $\partial M \neq \emptyset$.
Assume further  that $\pi : Y \to M$ has typical fibre $X$ which
is path connected, and put  $\pi_+ = \pi \amalg \id _M: Y \amalg M \to M$ be  
the $X_+$-bundle associated to $\pi$. Note that 
in this situation neither $(M, M_0)$ nor $X_+$
are connected. We  have
$$
C(M, M_0; \pi_+) = C(M; \pi_+) = \coprod_{k\geq 0} C_k (M; \pi),
$$
and though the scanning map is well-defined, it is no longer a (weak)  
homotopy equivalence. Similarly, the power set maps are still well-defined. But
as the filtration quotient $D^k$ is just $C_k$ the splitting theorems 
are trivially true. 
Nevertheless, we will now apply the power set maps  to  analyse the 
relation between the  spaces $C_k (M; \pi)$ as $k$ grows.

We  first define inclusion maps that allow us to think of the 
configurations of $k$ points as a subspace of the configuration 
space of $k+1$ points.
Recall $M \subset M^+ = M\cup \partial M\times [0,1/2)$ and $\pi $ 
can be extended to $M^+$ by 
pulling back along the natural projection $M^+ \to M$.
Fix points $z_0 \in M^+ \smallsetminus M$ and  $x_0 \in \pi ^{-1} (z_0)$.
We define a stabilisation map 
$$
b: C_k (M; \pi) \longrightarrow C_{k+1} (M^+; \pi)
	\longrightarrow C_{k+1} (M; \pi)
$$
where the first map adds the particle $(z_0,x_0)$ to any configuration 
and the second map is an injective homotopy equivalence  
induced by a map that isotopes $M^+$ into $M$ mapping 
$\partial M\times [0, 1/2)$ into a collar of $\partial M$ in $M$ while 
at the same time leaving the complement of the collar fixed.
As we assume that our diffeomorphisms always  fix a collar of $\partial M$ the 
following result follows immediately from these definitions.

\begin{proposition}
The inclusion $b: C_k (M;\pi) \to C_{k+1} (M; \pi)$ is $G$-equivariant. \qed
\end{proposition}

The fact that the map $b$ is stably split injective is well-known. 
The point of the next result is that this 
can be done $G$-equivariantly 
and that therefore we have stable splittings of Borel
constructions.

\begin{theorem}
The inclusion maps $b : C_k(M;\pi) \to C_{k+1} (M; \pi)$ are $G$-equivariantly,
stably split injective. In particular,
the  induced maps 
$$
EG \times _G C_k(M;\pi)\longrightarrow EG\times _G C_{k+1} (M; \pi)
$$
are stably split injective.
\end{theorem}

\begin{proof} 
Let $C^k = C_k (M; \pi)$. Then via the map  $b$ we consider $C^{k-1}$ to be a 
subspace of $C^k$ and $C = \lim _k C_k (M; \pi)$  
is a filtered space with
filtration quotient  $D^k = C^k/C^{k-1}$. Let  $V^k$ and $V$ be the 
wedge products of these. Since the inclusion map $b$ is $G$-equivariant, 
$G$ acts
naturally on $D^k$ and hence $V^k$ and $V$.

As in section 4.3 we can also define Borel configuration spaces $\mathscr C^k$,
$\mathscr C$, etc.  and a power 
set map 
$$
\sigma_k : \mathscr C^k = E\Sigma _k \times_{\Sigma _k} \widetilde 
C_k (M; \pi)  
\longrightarrow E \Sigma _{2^k} \times _{\Sigma_{ 2^k}} V^{2^k} 
$$ 
in this setting. As the action of $\Sigma _k$ on $\widetilde C_k$ is free, 
$\mathscr C^k$ is  ($G$-equivariantly)
homotopy equivalent to $C^k$ and hence $\mathscr C$ to $C$.
The arguments used in the proofs of Proposition 4.6 and 
Theorem 4.7 go through verbatim to show that the $\sigma _k$ fit together to
give a well-defined map of free $\Gamma ^+$-spaces
$$
\bar \sigma: \Gamma ^+ (\mathscr C) \longrightarrow \Gamma ^+ (V)
$$
which is $G$-equivariant, filtration preserving and a homotopy equivalence.
As both $\mathscr C$ and $V$ are connected, $\bar \sigma$ 
is homotopy equivalent 
to a map of free infinite loop spaces. This proves the fact that $\mathscr C$ 
and $V$ are stably homotopy equivalent.
 
Collapsing $D^k$ to a point defines a $G$-equivariant splitting of 
the inclusion $V^k \to V^{k+1}$. Via $\bar 
\sigma $ we conclude that  stably $\mathscr C^k \to \mathscr C ^{k+1} $ and
hence $b: C^k \to C^{k+1}$ are also  split injective.  
\end{proof}


Let $M \smallsetminus \bold k $ denote the manifold $M$ with a set of $k$
points deleted from the interior of $M$ (away from the collar of the boundary).
Using the same maps on the underlying manifolds as in the  
definition of $b$ we can define a homomorphism of diffeomorphism groups
$$
\bar b: \Diff (M \smallsetminus \bold{k}, \partial M) \longrightarrow
\Diff (M^+\smallsetminus \bold {k+1}, \partial M \times [0, 1/2) ) 
\longrightarrow  
\Diff (M \smallsetminus \bold {k+1}; \partial M).
$$

\begin{theorem}
The map $\bar b: B\Diff (M \smallsetminus \bold k, \partial M) 
\longrightarrow B\Diff (M \smallsetminus \bold {k+1}, \partial M )$
is stably split injective. In particular the induced map in homology is a
split injection in all degrees.
\qed
\end{theorem}

\begin{proof}
Fix $k$ distinct points ${\bold k} = \{ z_1, \dots z_k\}$ in the interior of 
$M$. 
Simultaneous
evaluation of diffeomorphisms on these points defines  a fibration
$$
\Diff (M \smallsetminus \bold k, \partial M )
\longrightarrow \Diff (M, \partial M) \longrightarrow C_k (M)
$$
and hence a fibration of classifying spaces
$$
C_k (M) \longrightarrow B\Diff (M \smallsetminus \bold k , \partial M )
\longrightarrow B\Diff (M, \partial M).
$$
Here we may take as  a model for the total space of the latter fibration  the 
Borel construction on
$C_k (M)$  by the following general construction: 
given a topological group $H$ with a closed  subgroup $K$
the fibration $K \to H \to H/K$ gives rise to a fibration 
$H/K \to EH \times _H H/K \simeq EH/ K \to EH/H =BH$  where $EH$ 
is a contractible space with a proper, free $H$-action; the total space of the latter fibration is a model for $BK$.  

With this identification, the result is now a special case of Theorem 4.14.
\end{proof}

\begin{remark}
When  $M$ is an orientable surface  $F _{g, n}^k$ of genus $g \geq 1$ 
with $n\geq 1$
boundary components and 
$k$ punctures this specialises to a result in \cite{BT}. 
The proof in \cite {BT} uses the
geometric construction of the power set map. Unfortunately, as defined there 
it is  not equivariant. 
The construction  needs to be replaced by the one in section 4.2 here. 
\end{remark}

In a forthcoming paper \cite{TIL} 
we will expand on the work here and 
show amongst other things that the maps $\bar b$ are furthermore
isomorphims
in homology in degrees $\leq k/2$ thus establishing a 
homology stability criteria for arbitrary manifolds and punctures. 



	\addcontentsline{toc}{section}{References}
	\bibliography{RM_Paper_Biblio}

\providecommand{\bysame}{\leavevmode\hbox to3em{\hrulefill}\thinspace}
\providecommand{\MR}{\relax\ifhmode\unskip\space\fi MR }
\providecommand{\MRhref}[2]{%
  \href{http://www.ams.org/mathscinet-getitem?mr=#1}{#2}
}
\providecommand{\href}[2]{#2}
\begin{thebibliography}{AAB80}

\bibitem[AAB80]{ABPM}
Ahmed~M. Abd-Allah and Ronald Brown, \emph{A compact-open topology on partial
  maps with open domain}, J. London Math. Soc. (2) \textbf{21} (1980), no.~3,
  480--486. 

\bibitem[BB78]{BBPM}
Peter~I. Booth and Ronald Brown, \emph{Spaces of partial maps, fibred mapping
  spaces and the compact-open topology}, General Topology and Appl. \textbf{8}
  (1978), no.~2, 181--195. 

\bibitem[BE74a]{BEG}
M.~G. Barratt and Peter~J. Eccles, \emph{{$\Gamma ^{+}$}-structures. {I}. {A}
  free group functor for stable homotopy theory}, Topology \textbf{13} (1974),
  25--45. 

\bibitem[BE74b]{BEG2}
\bysame, \emph{{$\Gamma ^{+}$}-structures. {II}. {A} recognition principle for
  infinite loop spaces}, Topology \textbf{13} (1974), 113--126. 


\bibitem[BM88]{BMHQ}
C.-F. B{\"o}digheimer and I.~Madsen, \emph{Homotopy quotients of mapping spaces
  and their stable splitting}, Quart. J. Math. Oxford \textbf{39} (1988),
  no.~156, 401--409. 

\bibitem[B{\"o}d87]{BSS}
C.-F. B{\"o}digheimer, \emph{Stable splittings of mapping spaces}, Algebraic
  topology ({S}eattle, {W}ash., 1985), Lecture Notes in Math., vol. 1286,
  Springer, Berlin, 1987, pp.~174--187. 

\bibitem[BT01]{BT}
Carl-Friedrich B{\"o}digheimer and Ulrike Tillmann, \emph{Stripping and
  splitting decorated mapping class groups}, Cohomological methods in homotopy
  theory ({B}ellaterra, 1998), Progr. Math., vol. 196, Birkh\"auser, Basel,
  2001, pp.~47--57.

\bibitem[DT58]{DTQ}
Albrecht Dold and Ren{{\'e}} Thom, \emph{Quasifaserungen und unendliche
  symmetrische {P}rodukte}, Ann. of Math. \textbf{67} (1958), 239--281.


\bibitem[Fra]{FRA}
John Francis, \emph{Factorization homology of topological manifolds}, Preprint;
  arXiv:1206.5522.

\bibitem[Gal11]{GAL}
S{\o}ren Galatius, \emph{Stable homology of automorphism groups of free
  groups}, Ann. of Math. (2) \textbf{173} (2011), no.~2, 705--768. 


\bibitem[God08]{GHST}
Veronique Godin, \emph{Higher string topology operations}, arXiv:0711.4859v2
  (2008).

\bibitem[GW99]{GW}
Thomas~G. Goodwillie and Michael Weiss, \emph{Embeddings from the point of view
  of immersion theory. {II}}, Geom. Topol. \textbf{3} (1999), 103--118
  (electronic). 

\bibitem[Hes92]{HHD}
Lars Hesselholt, \emph{A homotopy theoretical derivation of {$Q\;{\rm Map}\;
  (K,-)_+$}}, Math. Scand. \textbf{70} (1992), no.~2, 193--203. 
 

\bibitem[Lee03]{LSM}
John~M. Lee, \emph{Introduction to smooth manifolds}, Graduate Texts in
  Mathematics, vol. 218, Springer-Verlag, New York, 2003. 
  

\bibitem[Lur]{LUR}
Jacob Lurie, \emph{Higher algebra}, Preprint; available at
  http://www.math.harvard.edu/~lurie/.

\bibitem[May72]{MILS}
J.~P. May, \emph{The geometry of iterated loop spaces}, Springer-Verlag,
  Berlin, 1972, Lectures Notes in Mathematics, Vol. 271. 


\bibitem[McD75]{MDCS}
Dusa McDuff, \emph{Configuration spaces of positive and negative particles},
  Topology \textbf{14} (1975), 91--107. 

\bibitem[MW07]{MWA}
Ib~Madsen and Michael Weiss, \emph{The stable moduli space of {R}iemann
  surfaces: {M}umford's conjecture}, Ann. of Math. (2) \textbf{165} (2007),
  no.~3, 843--941.

\bibitem[Pal]{PAL}
Martin Palmer, \emph{Homology stability for configurations of submanifolds}, in
  DPhil-Thesis, Oxford 2013.

\bibitem[Sal01]{MR1851264}
Paolo Salvatore, \emph{Configuration spaces with summable labels},
  Cohomological methods in homotopy theory ({B}ellaterra, 1998), Progr. Math.,
  vol. 196, Birkh\"auser, Basel, 2001, pp.~375--395. 

\bibitem[Seg73]{SILS}
Graeme Segal, \emph{Configuration-spaces and iterated loop-spaces}, Invent.
  Math. \textbf{21} (1973), 213--221.

\bibitem[Sna74]{SSD}
V.~P. Snaith, \emph{A stable decomposition of {$\Omega ^{n}S^{n}X$}}, J. London
  Math. Soc. (2) \textbf{7} (1974), 577--583. 

\bibitem[Til]{TIL}
Ulrike Tillmann, \emph{Homology stability for symmetric diffeomorphisms and
  mapping classes}, in preparation.

\bibitem[Til12]{TILB}
\bysame, \emph{Spaces of graphs and surfaces: on the work of {S}\o ren
  {G}alatius}, Bull. Amer. Math. Soc. (N.S.) \textbf{49} (2012), no.~1, 73--90.
  

\end{thebibliography}
	\bibliographystyle{amsalpha}

\vskip .3in
Mathematical Institute
\newline
Oxford University
\newline
24-29 St Giles Street
\newline
Oxford OX1 3LB
\newline
UK

{\it manthorpe@maths.ox.ac.uk} and {\it tillmann@maths.ox.ac.uk}

\end{document}